\documentclass[11pt,reqno,final,fleqn]{amsart}\usepackage[]{graphicx}\usepackage{xcolor}
\makeatletter
\def\maxwidth{ %
  \ifdim\Gin@nat@width>\linewidth
    \linewidth
  \else
    \Gin@nat@width
  \fi
}
\makeatother

\definecolor{fgcolor}{rgb}{0.345, 0.345, 0.345}
\makeatletter
\@ifundefined{AddToHook}{}{\AddToHook{package/xcolor/after}{\definecolor{fgcolor}{rgb}{0.345, 0.345, 0.345}}}
\makeatother

\usepackage{framed}
\makeatletter
 {\par\unskip\endMakeFramed%
 \at@end@of@kframe}
\makeatother

\definecolor{shadecolor}{rgb}{.97, .97, .97}
\definecolor{messagecolor}{rgb}{0, 0, 0}
\definecolor{warningcolor}{rgb}{1, 0, 1}
\definecolor{errorcolor}{rgb}{1, 0, 0}
\makeatletter
\@ifundefined{AddToHook}{}{\AddToHook{package/xcolor/after}{
\definecolor{shadecolor}{rgb}{.97, .97, .97}
\definecolor{messagecolor}{rgb}{0, 0, 0}
\definecolor{warningcolor}{rgb}{1, 0, 1}
\definecolor{errorcolor}{rgb}{1, 0, 0}
}}
\makeatother
\newenvironment{knitrout}{}{} 

\usepackage{alltt}
\usepackage{tikz}

\newcommand\prob[1]{\mathbb{P}\left[{#1}\right]}
\newcommand\expect[1]{\mathbb{E}\left[{#1}\right]}
\newcommand\indicator[1]{\mathbbm{1}_{{#1}}}

\newcommand\lik{\mathcal{L}}

\newcommand{\dd}[1]{\mathrm{d}{#1}}
\newcommand{\halfopen}[1]{[{#1})}

\newcommand{\pd}[3][]{%
  \def\ord{#1} \ifx\ord\empty%
  \frac{\partial{#2}}{\partial{#3}}%
  \else \frac{\partial^{#1}{#2}}{\partial{#3}^{#1}}%
  \fi
}
\newcommand{\deriv}[3][]{%
  \def\ord{#1} \ifx\ord\empty%
  \frac{\dd{#2}}{\dd{#3}}%
  \else \frac{\dd^{#1}{#2}}{\dd{#3}^{#1}}%
  \fi
}
\newcommand{\opnorm}[1]{{\left\vert\kern-0.25ex\left\vert\kern-0.25ex\left\vert{#1}
        \right\vert\kern-0.25ex\right\vert\kern-0.25ex\right\vert}}

\newcommand\predec[1]{{#1}\mbox{-}}
\newcommand\prev[1]{{#1}\mbox{-}}
\newcommand\restrict[2]{{#1}\vert_{#2}}
\newcommand{\green}{\mathrm{green}}
\newcommand{\black}{\mathrm{black}}
\newcommand{\blue}{\mathrm{blue}}
\newcommand{\red}{\mathrm{red}}
\newcommand\Node[1]{\mathsf{#1}}
\newcommand\Ball[1]{\mathsf{#1}}

\def\Xspace{\mathbb{Z}^d}
\def\Uspace{\mathbb{Z}^d}
\def\N{\mathbb{N}}
\def\X{\mathcal{X}}
\def\Master{\mathcal{A}}
\def\Hist{\mathcal{H}}
\def\Gen{\mathcal{G}}
\def\Vis{\mathcal{V}}
\def\Wch{\mathcal{W}}
\def\Inv{\mathcal{I}}
\def\jumprate{\alpha}
\def\auxmeas{\beta}
\def\Borel{\mathcal{B}}
\def\Time{\mathbb{R}_{+}}
\def\rate{\mathbb{R}_{+}}
\def\cadlag{c\`adl\`ag\xspace}
\def\Bset{\mathbf{B}}
\def\Dset{\mathbf{D}}
\def\Sset{\mathbf{G}}

\def\Colors{\mathsf{F}}
\def\Eset{\mathsf{E}}
\def\Cset{\mathsf{C}}
\def\Lset{\mathsf{L}}
\def\Qset{\mathsf{D}}
\def\Masterset{\mathsf{A}}
\def\Tset{\mathsf{S}}
\def\Uset{\mathsf{U}}

\newcommand\seq[1]{\Node{P}({#1})}

\usepackage[round,elide]{natbib}
\usepackage{times}
\usepackage[neverdecrease]{paralist}
\usepackage{amssymb}
\usepackage{mathrsfs,eucal,bbm,mathtools}
\usepackage{xspace}
\usepackage[pdftex,colorlinks=true]{hyperref}

\usepackage{chngcntr}
\counterwithout{equation}{section}
\renewcommand\thesection{\arabic{section}}
\renewcommand\thesubsection{\thesection.\arabic{subsection}}

\usepackage[small,compact]{titlesec}
\newcommand\periodafter[1]{{#1}.}
\titleformat{\section}[hang]{\large\bfseries}{\periodafter\thesection}{2ex}{\periodafter}{}
\titleformat{\subsection}[hang]{\normalsize\bfseries}{\periodafter\thesubsection}{1ex}{\periodafter}{}
\titleformat{\subsubsection}[runin]{\normalsize\bfseries}{\thesubsubsection}{0em}{\periodafter}{}
\titleformat{\paragraph}[runin]{\normalsize\bfseries}{\theparagraph}{0em}{\periodafter}{}
\titlespacing{\section}{0em}{*1}{*0}
\titlespacing{\subsection}{0em}{*0}{*0}
\titlespacing{\paragraph}{0em}{*1}{*1}

\usepackage[sort&compress]{cleveref}

\theoremstyle{plain}
\newtheorem{thm}{Theorem}

\newtheorem{lemma}{Lemma}

\theoremstyle{definition}

\theoremstyle{remark}

\crefname{figure}{Fig.}{Figs.}
\Crefname{figure}{Fig.}{Figs.}
\crefname{table}{Table}{Tables}
\Crefname{table}{Table}{Tables}
\crefname{equation}{Eq.}{Eqs.}
\Crefname{equation}{Eq.}{Eqs.}
\crefname{pluralequation}{Eqs.}{Eqs.}
\Crefname{pluralequation}{Eqs.}{Eqs.}
\creflabelformat{equation}{#2#1#3}
\crefname{appendix}{Appendix}{Appendices}
\Crefname{appendix}{Appendix}{Appendices}
\crefname{chapter}{Ch.}{Chs.}
\Crefname{chapter}{Ch.}{Chs.}
\crefname{section}{\S}{\S\S}
\Crefname{section}{\S}{\S\S}
\crefformat{section}{#2\S#1#3}
\Crefformat{section}{#2\S#1#3}
\crefmultiformat{section}{#2\S\S#1#3}{ and~#2#1#3}{, #2#1#3}{, and~#2#1#3}
\Crefmultiformat{section}{#2\S\S#1#3}{ and~#2#1#3}{, #2#1#3}{, and~#2#1#3}
\crefname{thm}{Theorem}{Theorems}
\Crefname{thm}{Theorem}{Theorems}
\crefname{corol}{Corollary}{Corollaries}
\Crefname{corol}{Corollary}{Corollaries}
\crefname{prop}{Proposition}{Propositions}
\Crefname{prop}{Proposition}{Propositions}
\crefname{conj}{Conjecture}{Conjectures}
\Crefname{conj}{Conjecture}{Conjectures}
\crefname{lemma}{Lemma}{Lemmas}
\Crefname{lemma}{Lemma}{Lemmas}
\crefname{defn}{Definition}{Definitions}
\Crefname{defn}{Definition}{Definitions}
\crefname{inpar}{}{}
\Crefname{inpar}{Item}{Items}
\crefformat{inpar}{(#2#1#3)}

\usepackage{algorithm}
\usepackage{algorithmic}

\definecolor{darkgreen}{rgb}{0.0,0.5,0.2}

\definecolor{darkpurple}{rgb}{0.5,0.0,1.0}

\definecolor{darkbrown}{rgb}{0.8,0.8,0.0}

\title{Markov genealogy processes}
\author[King]{Aaron~A.~King}
\address{
  A.~A.~King,
  Department of Ecology \& Evolutionary Biology,
  Center for the Study of Complex Systems,
  Center for Computational Medicine \& Biology, and
  Michigan Institute for Data Science,
  University of Michigan,
  Ann Arbor, MI 48109 USA
}
\email{kingaa@umich.edu}
\urladdr{\href{https://kinglab.eeb.lsa.umich.edu/}{https://kinglab.eeb.lsa.umich.edu/}}
\author[Lin]{Qianying Lin}
\address{
  Q.-Y. Lin,
  Michigan Institute for Data Science,
  University of Michigan,
  Ann Arbor, MI 48109 USA
}
\author[Ionides]{Edward~L.~Ionides}
\address{
  E.~L.~Ionides,
  Department of Statistics and
  Michigan Institute for Data Science,
  University of Michigan,
  Ann Arbor, MI 48109 USA
}
\date{\today}

\hypersetup{pdftitle={Markov genealogy processes}}
\hypersetup{pdfauthor={A.A. King, Q.-Y. Lin, E.L. Ionides}}
\hypersetup{urlcolor=blue,citecolor=blue,linkcolor=blue,filecolor=blue}

\IfFileExists{upquote.sty}{\usepackage{upquote}}{}
\begin{document}

\maketitle

\begin{abstract}
  We construct a family of genealogy-valued Markov processes that are induced by a continuous-time Markov population process.
  We derive exact expressions for the likelihood of a given genealogy conditional on the history of the underlying population process.
  These lead to a nonlinear filtering equation which can be used to design efficient Monte Carlo inference algorithms.
  We demonstrate these calculations with several examples.
  Existing full-information approaches for phylodynamic inference are special cases of the theory.
\end{abstract}

\section{Introduction}
\label{sec:intro}

\emph{Phylodynamics} is the study of the traces left by the population dynamics of biological organisms in the genomes of their descendants, specifically, in the patterns of genealogical or phylogenetic relatedness among organisms or groups of organisms.
The term was first coined with reference to the processes of transmission and pathogen evolution in the context of infectious disease \citep{Grenfell2004,Frost2015} and to date many fruitful applications have been in this area \citep[e.g.,][]{Koelle2006a,ODea2011,Volz2013a,Rasmussen2014,Alizon2014,Faria2014,duPlessis2015,Geoghegan2015,Vijaykrishna2015,Biek2015,Smith2017,Li2017a,Hadfield2018,Bedford2020,RagonnetCronin2021}.
In typical applications, one wishes to infer the form and parameterization of a model of pathogen transmission on the basis of information contained in pathogen genome sequences.
Similar problems arise outside the realm of infectious disease biology, for example in systematics, comparative biology, cancer, microbiology, and population genetics \citep{Maddison2007,Gill2016,Stadler2021,MacPherson2021}.

A central technical problem in phylodynamics is the precise characterization of the relationship between stochastic population processes and the traces they leave in genomes obtained from a sample of the population.
This in turn is usually factored into two subproblems:
\begin{inparaenum}[(i)]
\item the relationship between genome sequences and the genealogies or phylogenies that relate them and
\item the relationship between these genealogies and the population processes that generate them.
\end{inparaenum}
In this paper, we focus attention on the second subproblem.
In particular, we suppose that we have one or more stochastic models that describe the dynamics of infections and recoveries (or births and deaths, or speciations and extinctions) in a population
and we desire to estimate parameters and compare these models using data in the form of genealogies that relate sampled infections, individuals, or species.

At present, three broad approaches to the solution of this problem exist.
The first compares the reconstructed genealogy with simulated genealogies on the basis of summary statistics designed to capture the important features of the genealogy.
Such approaches can be used to estimate parameters and compare models, for example via approximate Bayesian computation \citep{Sisson2007,Luciani2009,Ratmann2012,Poon2015} or synthetic likelihood \citep{Wood2010,Fasiolo2016}.
These methods have the advantage of being applicable for essentially arbitrary models, but cast away some of the information contained in the data.
It can be difficult to quantify the amount of information discarded and to design summary statistics that minimize this information loss.
By contrast, the second and third approaches are \emph{full-information} methods, in the sense that they are based on the likelihood.
The approach pioneered by Volz, Koelle, and colleagues \citep{Volz2009a,Rasmussen2011,Koelle2012,Volz2013,Dearlove2013} is based on the \citet{Kingman1982,Kingman1982c,Kingman1982b} coalescent, which yields the distribution of genealogies as a function of a time-varying \emph{coalescent rate}.
These approaches rely on the assumption that the population size is large and the sample size small, so that branching points in the sample lineages are approximately independent of those in the population.
Rigorous quantification of the error associated with this approximation exists only for special cases \citep{Fu2006}.
The third approach, associated with Stadler and colleagues \citep{Gernhard2008,Stadler2010,Stadler2012,Boskova2014,Kuehnert2014,MacPherson2021}, is based on birth-death processes.
In the case of the linear birth-death process, exact expressions can be obtained, but approximations must be employed to deal with model nonlinearity.
In particular, reverse-time arguments that go through in the linear case fail in the nonlinear case.

Recently, \citet{Etheridge2019} described an approach similar in objective to the one described here, but based instead on Fleming-Viot processes and lookdown constructions \citep{Donnelly1996,Donnelly1999}.
More closely related to our approach is the work of \citet{Vaughan2019}, who recently devised an algorithm for exact phylodynamic likelihood computation for simple jump processes under certain assumptions of time-homogeneity.
The arguments in this paper put this algorithm on a firm footing and place it in a much broader context, allowing consideration of a wider class of models and laying the groundwork for more efficient algorithms.
We return to this issue in \cref{sec:discussion}.

In this paper, we define the notion of the \emph{genealogy process} that is induced by a population model, i.e., a model of the births and deaths occurring in a specified population.
In the general case, these births and deaths will be stochastic, which implies that the induced genealogy process will itself be a stochastic process on the space of genealogies.
Accordingly, we are interested in characterizing the probabilistic properties of this process.
In practice, information about any real population is obtained through genomic sampling, and we will therefore also be interested in the properties of partial genealogies that relate the samples.

As for the population models, we will restrict our attention to Markov processes.
In practice, this is not a major restriction, as most of the models of scientific interest can be readily formulated as, or approximated by, Markov processes.
To make the notion of a genealogy definite, it is necessary to conceive of births and deaths as discrete events.
It is natural, therefore, to further restrict our attention to Markov jump processes, i.e., continuous-time processes for which births and deaths occur at random times.
The approach we describe here can be generalized to a somewhat broader category of Markov processes, but the novel mathematical constructions are already fairly complex, so we postpone these generalizations to a future contribution.
In particular, we simplify matters by confining ourselves to the case where birth, death, and sampling events occur one at a time almost surely.
We note that some models of theoretical and practical interest do violate this assumption and that our approach can be generalized to accommodate such models.
Nevertheless, we defer consideration of these interesting complexities to a sequel.

The remainder of the paper is structured as follows.
In \cref{sec:prelims}, we lay mathematical groundwork by constructing a probability space within which we can speak of the genealogies induced by a given population process.
In \cref{sec:examples}, we give several examples, both of processes amenable to treatment within this framework, and of models that motivate further theoretical development.
Next, in \cref{sec:genealogies}, we define several related Markov genealogy processes induced by a given population model.
With these definitions in hand, in \cref{sec:results} we derive our main results.
There, we exhibit explicit expressions for the likelihood of an observed genealogy and derive the analogue of the Duncan-Mortensen-Zakai (DMZ) or nonlinear filtering equation for these processes.
In \cref{sec:illus}, we revisit some of the examples of \cref{sec:examples} to demonstrate these likelihood computations more concretely.
Finally, in \cref{sec:discussion}, we indicate some of the implications for phylodynamics broadly.
In particular, we note that both major existing approaches for likelihood-based phylodynamic inference are special cases.
We also point out that, from the DMZ equation, a straight road leads to efficient Monte Carlo inference algorithms capable of accommodating a broader class of models than has heretofore been susceptible to analysis.

\section{Mathematical preliminaries}
\label{sec:prelims}

\paragraph{Markov jump processes on the integer lattice}

Suppose we have a non-explosive Markov jump process $\X_t\in\Xspace$, parameterized by $t\in\Time$, which we think of as a model of the random time-evolution of some kind of population.
Henceforth, we refer to $\X_t$ as a \emph{population process}.
It is defined by its initial-state distribution and its generator.
In particular, we suppose that
\begin{equation}\label{eq:ic}
  \prob{\X_0=x}=p_0(x),
\end{equation}
for some choice of initial distribution, $p_0$.
The transitions of $\X_t$ are governed by event-rate functions $\jumprate_u(t,x)\in\rate$, for $u,x\in\Xspace$, $t\in\Time$.
In particular, $\jumprate_u(t,x)$ is the hazard of a jump from state $x$ to state $x+u$ at time $t$.
These conditions imply that, for any $f\in L^{\infty}(\Xspace)$, if $F(s,x)\coloneqq\expect{f(\X_t)\;|\;\X_s=x}$ for $0\le{s}\le{t}$, then $F$ satisfies the final-value problem
\begin{equation}\label{eq:gen}
  \pd{F}{s}(s,x)=-\sum_{u\in\Xspace}\!\jumprate_u(s,x)\,\left[F(s,x+u)-F(s,x)\right],\qquad F(t,x) = f(x), \qquad x\in\Xspace.
\end{equation}
If we moreover assume that the sample paths of $\X_t$ are right-continuous with left limits (\cadlag), then \cref{eq:ic,eq:gen} completely specify $\X_t$.
The assumption that $\X_t$ is non-explosive, and the further requirement we make that $\sum_u\jumprate_u(t,x)<\infty$ for all $t$ and $x$ restricts the class of allowable rate functions $\jumprate$.

The adjoint of \cref{eq:gen} is the Kolmogorov forward equation (sometimes called the \emph{master equation} in this context):
\begin{equation}\label[pluralequation]{eq:kfe}
  \begin{gathered}
    \pd{w}{t}(t,x)=\sum_{u\in\Xspace}\!\jumprate_u(t,x-u)\,w(t,x-u)-\jumprate_u(t,x)\,w(t,x),\qquad x\in\Xspace,\\
    w(0,x) = p_0(x).
  \end{gathered}
\end{equation}
If $w(t,x)$ satisfies \cref{eq:kfe}, then $w(t,x)=\prob{\X_t=x}$.

\paragraph{Definitions}

\begin{figure}
  \resizebox{0.6\linewidth}{!}{
    \begin{tikzpicture}[scale=1.4]
      \usetikzlibrary{arrows}
      \tikzstyle{coordinate}=[inner sep=0pt,outer sep=0pt]
      \tikzstyle{det}=[color=black, thick, >=stealth]
      \tikzstyle{stoch}=[color=black, dashed, thick, >=stealth]
      \node (O) at (0,0) {$\Master_t$};
      \node (X) at (180:2) {$\X_t$};
      \node (H) at (216:2) {$\Hist_t$};
      \node (I) at (252:2) {$\Inv_t$};
      \node (G) at (288:2) {$\Gen_t$};
      \node (V) at (324:2) {$\Vis_t$};
      \node (W) at (360:2) {$\Wch_i$};
      \draw [det,->] (O) -- (X);
      \draw [det,->] (O) -- (H);
      \draw [det,->] (H) -- (X);
      \draw [det,->] (O) -- (I);
      \draw [det,->] (O) -- (G);
      \draw [det,->] (O) -- (V);
      \draw [det,->] (O) -- (W);
      \draw [det,->] (V) -- (W);
      \draw [det,->] (G) -- (I) node[midway,below,sloped] {$\mathrm{ext}$};
      \draw [det,->] (G) -- (V) node[midway,below,sloped] {$\mathrm{obs}$};
      \draw [stoch,->] (H) -- (I);
      \draw [stoch,->] (H) -- (G);
      \draw [stoch,->] (H) -- (V);
      \draw [stoch,->] (H) -- (W);
    \end{tikzpicture}
  }
  \caption{
    \label{fig:constellation}
    Relations among the various Markov processes discussed in the paper.
    Deterministic maps are indicated with solid arrows;
    random maps are shown as dashed arrows.
    All the maps shown commute.
    $\X_t$ is the \emph{population process}, a model of the dynamics of some system, which we take as a starting point.
    $\Hist_t$ is the \emph{history process}, which records the full history of $\X_t$.
    $\Inv_t$ is the \emph{inventory process}: at each time $t$, $\Inv_t$ is an inventory of all extant individuals in the population, each of which has a globally unique name.
    $\Gen_t$ is the \emph{genealogy process}, which captures the precise genealogical relationships among all individuals in $\Inv_t$, as well as among any samples that have been taken from the population.
    $\Vis_t$ is the \emph{visible genealogy process}, which is $\Gen_t$ pruned so that only relationships among samples remain.
    Finally $\Wch_i$ is the \emph{embedded chain of the visible genealogy process}, which is $\Vis_{s_i}$, $s_i$ being the time of the $i$-th sample.
    All of these processes can be obtained via deterministic procedures applied to the \emph{master process} $\Master_t$, as described in the text.
  }
\end{figure}
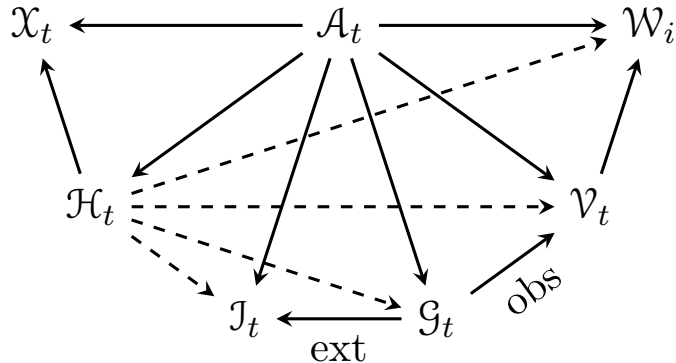

Our goal in this paper is to introduce a family of Markov processes induced by the jump process just described (\cref{fig:constellation}).
While population processes of the kind described above can be constructed in a variety of time-honored ways \citep[e.g.,][]{Andersen1993,Kallenberg1997}, these classical approaches are not entirely sufficient for the tree-valued Markov processes we will erect on top of the population process.
We therefore explicitly construct the probability space that underlies the stochastic processes we subsequently describe.
Unavoidably, this leads to technicalities that are necessary for the firm establishment of the properties of these processes but that may distract from the overarching goals.
Readers willing to stipulate these properties may skim the remainder of this section, in which we formally construct the probability space and make some definitions that will be needed in the sequel.

Let us define a \emph{jump} to be an ordered triple $(t,u,n)\in\Time\times\Uspace\times\N$.
We refer to $t$ as the \emph{time} of the jump;
$u$ is the \emph{type} of the jump;
$n$ is an \emph{auxiliary number} whose use will be made clear below.
A \emph{jump sequence} is a countable sequence of jumps at increasing times.
That is,
\begin{equation*}
  \omega=\left(t_k,u_k,n_k\right)_{k=0}^{K}
\end{equation*}
is a jump sequence if and only if $K\in\N\cup\{\infty\}$, $u_k\in\Uspace$, $n_k\in\N$ for all $k$, and $0=t_0<t_1<t_2<\dots$.
We will take our sample space, $\Omega$, to be the set of all jump sequences.
For $\omega\in\Omega$ as above, we write
\begin{equation}\label{eq:components1}
  \begin{gathered}
    T_k(\omega)\coloneqq t_k, \qquad
    U_k(\omega)\coloneqq u_k, \qquad
    N_k(\omega)\coloneqq n_k, \qquad
    K(\omega)\coloneqq K
  \end{gathered}
\end{equation}
and also
\begin{equation}\label{eq:components2}
  \begin{gathered}
    T(\omega)\coloneqq\left(T_k(\omega)\right)_{k=0}^{K(\omega)}, \qquad
    U(\omega)\coloneqq\left(U_k(\omega)\right)_{k=0}^{K(\omega)}, \qquad
    N(\omega)\coloneqq\left(N_k(\omega)\right)_{k=0}^{K(\omega)}.
  \end{gathered}
\end{equation}
We will denote by $\mathring{\Omega}$ the set of all finite jump sequences, i.e., $\mathring{\Omega}\coloneqq\{\omega\in\Omega\;|\;K(\omega)<\infty\}$.

Note that every element of $\Omega$ corresponds to a unique sample path of $\X_t$.
In particular, with the convention that $U_0=\X_0$, we will write
\begin{equation}\label{eq:Xdef}
  \X_t(\omega) = \sum_{k=0}^{K(\omega)}\!U_k(\omega)\,\indicator{\halfopen{T_k(\omega),\infty}}(t).
\end{equation}

The sample space $\Omega$ has a natural partial order.
We write $\omega\preceq\omega'$ if $\omega'$ is an extension of $\omega$;
that is, if $K(\omega)\le K(\omega')$, and $\left(T_k(\omega),U_k(\omega),N_k(\omega)\right)=\left(T_k(\omega'),U_k(\omega'),N_k(\omega')\right)$ for $k=0,\dots,K(\omega)$.
For $\omega\in\Omega$, the set
$\left\{\omega'\in\Omega\;\vert\;\omega'\preceq\omega\right\}$
is totally ordered.
Moreover, for each $\omega\in\mathring{\Omega}$, there is a unique \emph{predecessor}, $\predec{\omega}$, such that $\predec{\omega}\prec\omega$ and, for all $\omega'$, $\predec{\omega}\preceq\omega'\preceq\omega$ implies that either $\omega'=\predec{\omega}$ or $\omega'=\omega$.

In order to define probabilistic \emph{events}, it is necessary to define the $\sigma$-algebra of measurable subsets of $\Omega$.
We topologize $\Omega$ using a standard approach that makes $\Omega$ separable, i.e., possessed of a countable dense subset.
Specifically, let $d_S$ represent the Skorokhod metric on the space of \cadlag functions $\Time\to\Xspace$ \citep{Kallenberg1997,Ethier2009}.
Extend this to a metric $d$ on $\Omega$ by
\begin{equation*}
  d(\omega,\omega')\coloneqq d_S(\X(\omega),\X(\omega'))+\sup_{k}\left\vert N_k(\omega)-N_k(\omega')\right\vert,
\end{equation*}
where the sample-path functions $\X_t$ are defined above (\cref{eq:Xdef}) and it is understood that $N_k(\omega)=0$ for $k>K(\omega)$.
It is straightforward to verify that this is indeed a metric and that, equipped with this metric, $\Omega$ is a complete, separable metric space, i.e., a so-called Polish space.
We take our event space, $\Borel$, to be the Borel $\sigma$-algebra of $\Omega$.

For $t\in\Time$ and $\omega\in\Omega$, we define the \emph{restriction}, $\restrict{\omega}{t}\in\Omega$, by
\begin{equation*}
  \restrict{\omega}{t} \coloneqq \max\left\{\omega'\in\Omega\;\vert\;\omega'\preceq\omega\ \text{and}\ T_{K(\omega')}\le t\right\}.
\end{equation*}
and let $\restrict{\Omega}{t}$ be the space of $t$-restrictions.
Note that each $\restrict{\Omega}{t}$ is a complete, separable metric space.

We now turn to the probability measure, $\mathbb{P}$, on $\Omega\cong T(\Omega){\times}U(\Omega){\times}N(\Omega)$.
We will specify this by giving its density with respect to a base measure.
A natural base measure on $T(\Omega)$ is the probability measure of the rate-$\mu$ Poisson point process \citep{Andersen1993,Kallenberg1997};
the counting measure serves for the discrete component $U(\Omega){\times}N(\Omega)$.
Let $\pi_{\mu}$ denote the product of these two measures.
We define a probability measure on $\Omega$ by specifying its density with respect to $\pi_{\mu}$.
In particular, for each $u,x\in\Xspace$, we suppose $\auxmeas_{u,x}$ is a given probability measure on $\N$;
these will take specific forms below.
For $t\in\Time$ and $\omega\in\restrict{\Omega}{t}$, define the probability density function
\begin{equation}\label{eq:probmeas}
  \begin{aligned}
    P_t(\omega)\:\coloneqq\:&p_0(\X_0)\,\prod_{k=1}^{K}\left(\frac{1}{\mu}\,\jumprate_{U_k}(T_k,\X_{T_k}-U_k)\,\auxmeas_{U_k,\X_{T_k}-U_k}(N_k)\right)\\
    &\times\exp{\left(\mu\,t-\int_0^t\sum_{u\in\Xspace}{\jumprate_u(s,\X_s)}\,\dd{s}\right)}.
  \end{aligned}
\end{equation}
\normalsize
Here, for the sake of readability, we have suppressed the dependence of the random variables $\X_t$, $K$, $T_k$, $U_k$ and $N_k$ on $\omega$.
It is readily verified that, for all $t$, $P_t$ is a probability density with respect to $\pi_{\mu}$.
One can dispense with the dependence of \cref{eq:probmeas} on the Poisson rate parameter $\mu$.
For example, one can always choose $\mu=1$
and the factor of $e^t$, that remains can usually be neglected, it being a mere normalizing constant.
However, we preserve the dependence on $\mu$ here to remind us of the manner in which the magnitude of $P_t$ depends on the base measure, $\pi_{\mu}$.

\paragraph{Master process}

We now define a process that serves as the foundation for all that follows.
For $t\in\Time$ and $\omega\in\Omega$, let
\begin{equation}\label{eq:Adef}
  \Master_t(\omega) \coloneqq \left(t, \restrict{\omega}{t}\right).
\end{equation}
We will refer to $\Master_t$ as the \emph{master process}.
We use the expressions $t(\Master_t)$ and $\omega(\Master_t)$ to refer to the first and second elements of $\Master_t$, respectively.

We define the densities, $P_{\Master_{t_1},\dots,\Master_{t_m}}$, of finite collections of random variables $\{\Master_{t_1},\dots,\Master_{t_m}\}$, $t_1<t_2<\dots<t_m$, in terms of \cref{eq:probmeas} by defining,
\begin{equation}\label{eq:fd_dens}
  P_{\Master_{t_1},\dots,\Master_{t_m}}(a_1,\dots,a_m)\coloneqq
  \begin{cases}
    P_{t_m}(\omega(a_m)), &\text{if}\ \omega(a_1)\preceq\omega(a_2)\preceq\cdots\preceq\omega(a_m),\\
    0, &\text{otherwise}.
  \end{cases}
\end{equation}
Since the densities of \cref{eq:fd_dens} are projective and the state-space is Polish, the Kolmogorov Extension Theorem implies that there is a unique probability measure $\mathbb{P}$ on $\Omega$ with these finite-dimensional densities.
Thus $(\Omega,\Borel,\mathbb{P})$ is a probability space.
With these definitions, one readily verifies that the master process $\Master_t$ is Markov and that the population process $\X_t$, defined by \cref{eq:Xdef}, coincides with the Markov jump process defined by \cref{eq:ic,eq:gen}.
Moreover, the non-explosion assumption guarantees that $\omega(\Master_t)\in\mathring\Omega$ for every $t$.


\paragraph{History process}

Next among our constellation of related processes (\cref{fig:constellation}) is the \emph{history process}, which encapsulates the entire history of $\X$ up to time $t$.
Specifically, for $\omega\in\Omega$, we define
\begin{equation*}
  \Hist_t(\omega) \coloneqq \left(t,\left(T_k(\omega),U_k(\omega)\right)_{k=0}^{K(\restrict{\omega}{t})}\right),
\end{equation*}
where $T_k$, $U_k$, and $K$ are as in \cref{eq:components1}.
Thus $\Hist_t$ contains exactly those elements of $\omega$ that are relevant to the construction of $\X_t$.
It is trivial to verify that $\Hist_t$ is Markov and to compute its probability density.
In particular, given any history $h_t = \left(t,\left(t_j,u_j\right)_{j=0}^{k}\right)$,
the marginal density at $h_t$ is obtained by summing \cref{eq:probmeas} over all possible values of the finite sequence $N(\restrict{\omega}{t})$, which yields
\begin{equation*}
  P_{\Hist_t}(h_t)\coloneqq p_0(x_0)\,\prod_{j=1}^{k}\left(\frac{\jumprate_{u_j}(t_j,x_{t_j}-u_j)}{\mu}\right)\,\exp{\left(\mu\,t-\int_0^t\sum_{u\in\Xspace}{\jumprate_u(s,x_s)}\,\dd{s}\right)},
\end{equation*}
where, according to \cref{eq:Xdef}, $x_s=\sum_{j=0}^{k}u_j\,\indicator{\halfopen{t_j,\infty}}(s)$.
Then the probability density of $\Master_t$ conditional on $\Hist_t=h_t$ is
\begin{equation}\label{eq:condH}
  P_{\Master_t|\Hist_t}(\Master|h_t)\coloneqq\prod_{j=1}^{k}\auxmeas_{u_j,x_{t_j}-u_j}(N_j(\omega(\Master))).
\end{equation}
That is, conditional on the history process, the auxiliary numbers $N_j$ are independent random variables.

\paragraph{Births, deaths, population size}

The population process, $\X_t$ we have defined will, in applications, track the time-evolution of a structured population composed of discrete, exchangeable individuals.
In particular, we are interested in the genealogical relationships among members of some focal subpopulation.
For example, if we are interested in viral pathogen genealogies, the focal subpopulation might be the population of infected hosts.
In the event we are studying species phylogenies, the focal subpopulation might be a group of related species.
To facilitate this, we give some additional structure to the probability space we have constructed.
In particular, we will suppose there are functions $I,B,D:\Xspace\to\N$ such that 
\begin{equation}\label{eq:BDcondn}
  \jumprate_u(t,x)>0 \implies  I(x+u)-I(x)=B(u)-D(u),
\end{equation}
for all $x,u\in\Xspace$.
For any $x\in\Xspace$, $I(x)$ represents the size of the focal subpopulation when $\X_t=x$.
We interpret $B(u)$ as the number of births into our focal subpopulation associated with an event of type $u$, and $D(u)$ as the number of deaths.
\Cref{eq:BDcondn} guarantees that $I$ is compatible with $B$ and $D$:
it implies that the difference in population sizes between state $x$ and state $y$ matches the sum of births minus deaths over any possible path from $x$ to $y$.

Although it is both interesting and possible to treat the general case, in this paper, we assume that births and deaths occur one at a time and never co-occur.
In particular, we assume that $B(\Xspace),D(\Xspace)\subseteq\{0,1\}$.
Let $\Bset\coloneqq B^{-1}(\{1\})$ and $\Dset\coloneqq D^{-1}(\{1\})$ be the sets of event-types associated with births and deaths, respectively.
Our insistence that births and deaths not co-occur implies that $\Bset\cap\Dset=\emptyset$.

\paragraph{Samples}

Similarly, we suppose there is a function $G:\Xspace\to\{0,1\}$ that we interpret as the number of samples associated with each type of event.
That is, $u\in\Sset\coloneqq G^{-1}(\{1\})$ implies that an event of type $u$ results in a sample being taken.
As with births and deaths, we suppose $\Sset\cap(\Bset\cup\Dset)=\emptyset$.
Again, it is possible to relax the assumptions that samples occur singly and do not coincide with births or deaths, but the resulting technical complexities are best handled after the simpler case is in hand.

\paragraph{Absence of structure within the focal subpopulation}

Note that we have assumed the existence of a single focal subpopulation.
We will need the additional assumption that this focal subpopulation is itself \emph{unstructured}.
In particular, we will suppose that the individuals in the focal subpopulation (of size $I(\X)$) are exchangeable, in the sense that each is as likely as any other of being parent to the next newborn, of being sampled, or of dying.
To enforce this assumption, we postulate that, for $u\in\Bset\cup\Dset\cup\Sset$, $\auxmeas_{u,x}$ is uniform on the first $I(x)$ natural numbers, i.e.,
\begin{equation*}
  \auxmeas_{u,x}(n)=\frac{\indicator{\halfopen{0,I(x)}}(n)}{I(x)}.
\end{equation*}

\section{Examples}
\label{sec:examples}

Many interesting models fit within the constraints we have so far described.
Here, we enumerate a few familiar ones.
We also provide some examples of models that do not conform to our assumptions.

\paragraph{SIR and SIRS models}

The most basic model of an immunizing infection is the SIR model, whereby susceptible individuals immediately become infectious upon being infected and remain so until they recover or are otherwise removed from the population and recovered individuals are permanently immunized against reinfection.
Relaxing the latter assumption, i.e., allowing for waning of immunity, leads to the SIRS model, a modest extension.
For these models, we take $d=4$, so that the state vector is $\X=(s,i,r,g)$, $s$ representing the number of susceptibles in the population;
$i$, the number of infectives;
$r$, the number of recovered and immune hosts;
and $g$ the cumulative number of genomic samples collected.
There are four kinds of jumps, so that the rate function is
\begin{equation*}
  \alpha_u=
  \begin{cases}
    b(t)\,s\,i, & u=(-1,1,0,0),\ s>0,\ i>0,\\
    \gamma\,i,  & u=(0,-1,0,0),\ i>0,\\
    \sigma\,r,  & u=(1,0,-1,0),\ r>0,\\
    \psi(t,s,i,g)\,i, & u=(0,0,0,1),\ i>0,\\
    0, & \text{otherwise}.\\
  \end{cases}
\end{equation*}
Here, the rates shown are those of infection, recovery, loss of immunity, and sampling, respectively.
The parameters $\gamma$ and $\sigma$ are in this case constants, but we allow for a time-varying transmission rate, $b(t)$.
The sampling rate $\psi$ is an arbitrary function, except for the minimal constraints on the rate functions described in \cref{sec:prelims}.
In this case, the population size function is simply $I(\X)=i$, whilst $\Bset=\{(-1,1,0,0)\}$, $\Dset=\{(0,-1,0,0)\}$, and $\Sset=\{(0,0,0,1)\}$.

\paragraph{S\textsuperscript{2}IR model}

There is no barrier to having more than one susceptible class representing, for example, different risk groups.
For example, if we have two different susceptible classes, such that the \emph{per capita} risk of infection in the first is $b_1\,i$ and that in the second is $b_2\,i$, $i$ being the number of infectious individuals, then we have $d=4$, $\X=(s_1,s_2,i,g)$, where $s_1$, $s_2$ are the numbers in each of the two susceptible classes and $g$ is as before.
The jump rates are
\begin{equation*}
  \alpha_u=
  \begin{cases}
    b_1(t)\,s_1\,i, & u=(-1,0,1,0),\ s_1>0,\ i>0\\
    b_2(t)\,s_2\,i, & u=(0,-1,1,0),\ s_2>0,\ i>0\\
    \gamma\,i, & u=(0,0,-1,0),\ i>0\\
    \psi(t,s_1,s_2,i,g)\,i, & u=(0,0,0,1),\ i>0\\
    0, & \text{otherwise}.
  \end{cases}
\end{equation*}
The four non-zero rates are those of infection from the first class, infection from the second class, recovery, and sampling, respectively.
In this case, the population size function is again $I(\X)=i$, whilst $\Bset=\{(-1,0,1,0),(0,-1,1,0)\}$, $\Dset=\{(0,0,-1,0)\}$, and $\Sset=\{(0,0,0,1)\}$.

\paragraph{Linear birth-death-sampling process}

Linear processes have proved useful as models of speciation and extinction \citep{Nee1994a,Maddison2007,Gernhard2008,Tavare2018}.
For example, if we take $d=2$, $\X=(n,g)$, $I(\X)=n$, and
\begin{equation*}
  \alpha_u=
  \begin{cases}
    \lambda(t)\,n, & u=(1,0),\ n>0,\\
    \delta(t)\,n, & u=(-1,0),\ n>0,\\
    \psi(t)\,n, & u=(0,1),\ n>0,\\
    0, & \text{otherwise},
  \end{cases}
\end{equation*}
where $\lambda$, $\delta$, $\psi$ are the per-capita birth, death, and sampling rates, respectively.
If these are constants, then
we obtain the linear birth-death process with a constant \emph{per capita} sampling rate treated by \citet{Stadler2010}.
Of course, our formalism embraces nonlinear birth-death-sampling processes as well.

\paragraph{SEIR and SI\textsuperscript{2}R model}

As mentioned above, we will rely on the assumption that the individuals whose genealogies we study are exchangeable, in the sense that each is as likely as any other of giving birth, being sampled, or of dying.
This precludes consideration of a number of interesting models, including models with a latent period and models with heterogeneity of infectiousness.

\paragraph{Moran process}

We have explicitly ruled out the possibility that birth and death events co-occur, which prevents us from applying the theory we develop here to the classical Moran model which plays such an important role in population genetics \citep{Moran1958,Wakeley2008}.
We do so only to avoid distracting technicalities in this initial presentation.
In fact, with minor modifications, the theory can be extended to deal with this case as well as situations where sampling events co-occur with death events \citep[as in][]{Leventhal2014}.
We postpone consideration of these cases to a later paper.

\paragraph{Superspreading events}

Likewise, we defer consideration of jump processes for which multiple birth or death events can occur simultaneously.
Such processes deserve consideration in their own right (for example, as models of superspreading events) and as models of overdispersed population processes \citep{Breto2011}.
Again, accommodating such processes requires only a modest extension of the present theory, but the notational complexity introduced thereby recommends postponement of these developments to a forthcoming paper.
Thus, the genealogical trees under consideration in this paper will necessarily be binary.

\section{Markov genealogy processes}
\label{sec:genealogies}

\paragraph{Inventory process}

It will be helpful to introduce another Markov process that tracks the composition of the population through time.
We assume that every individual born into our focal subpopulation has a unique name (or more accurately, serial number), so that at every instant, the composition of the population is characterized by an \emph{inventory}, i.e., a list of the names of all extant (i.e., currently living) individuals.
When a birth occurs, the inventory is augmented with the new, globally unique name;
when a death occurs, one name is struck off the list.

To be precise, we define $\Inv_t(\omega)\coloneqq\mathrm{inven}(\restrict{\omega}{t})$, where $\mathrm{inven}$ is a deterministic, recursive procedure as follows.
First, define the $\mathrm{add}$ and $\mathrm{drop}$ operations:
if $\Inv=\{c_0,\dots,c_{m-1}\}\subset\N$, with $c_0<\dots<c_{m-1}$ and $n<m$,
then $\mathrm{add}(\Inv)=\Inv\cup\{c_{m-1}+1\}$ and $\mathrm{drop}(\Inv,n)=\Inv\setminus\{c_n\}$.
With these definitions, we write, for $\omega\in\mathring\Omega$,
\begin{equation}\label{eq:Idef}
  \mathrm{inven}(\omega)\coloneqq
  \begin{cases}
    \left\{0,\dots,I(\X_0(\omega))-1\right\}, &\text{if}\ K(\omega)=0;\\
    \mathrm{add}(\mathrm{inven}(\predec{\omega})), &\text{if}\ K(\omega)>0\ \text{and}\ U_{K(\omega)}\in\Bset;\\
    \mathrm{drop}(\mathrm{inven}(\predec{\omega}),N_{K(\omega)}(\omega)), &\text{if}\ K(\omega)>0\ \text{and}\ U_{K(\omega)}\in\Dset;\\
    \mathrm{inven}(\predec{\omega}), &\text{otherwise}.
  \end{cases}
\end{equation}
In view of \cref{eq:BDcondn}, it is clear that, for all $t$, $|\Inv_t|=I(\X_t)$.
Note that the non-explosiveness assumption is needed to guarantee the existence of $\predec{\omega}$.

\begin{figure}
\begin{knitrout}\small
\definecolor{shadecolor}{rgb}{0.969, 0.969, 0.969}\color{fgcolor}

{\centering \includegraphics[width=0.7\linewidth]{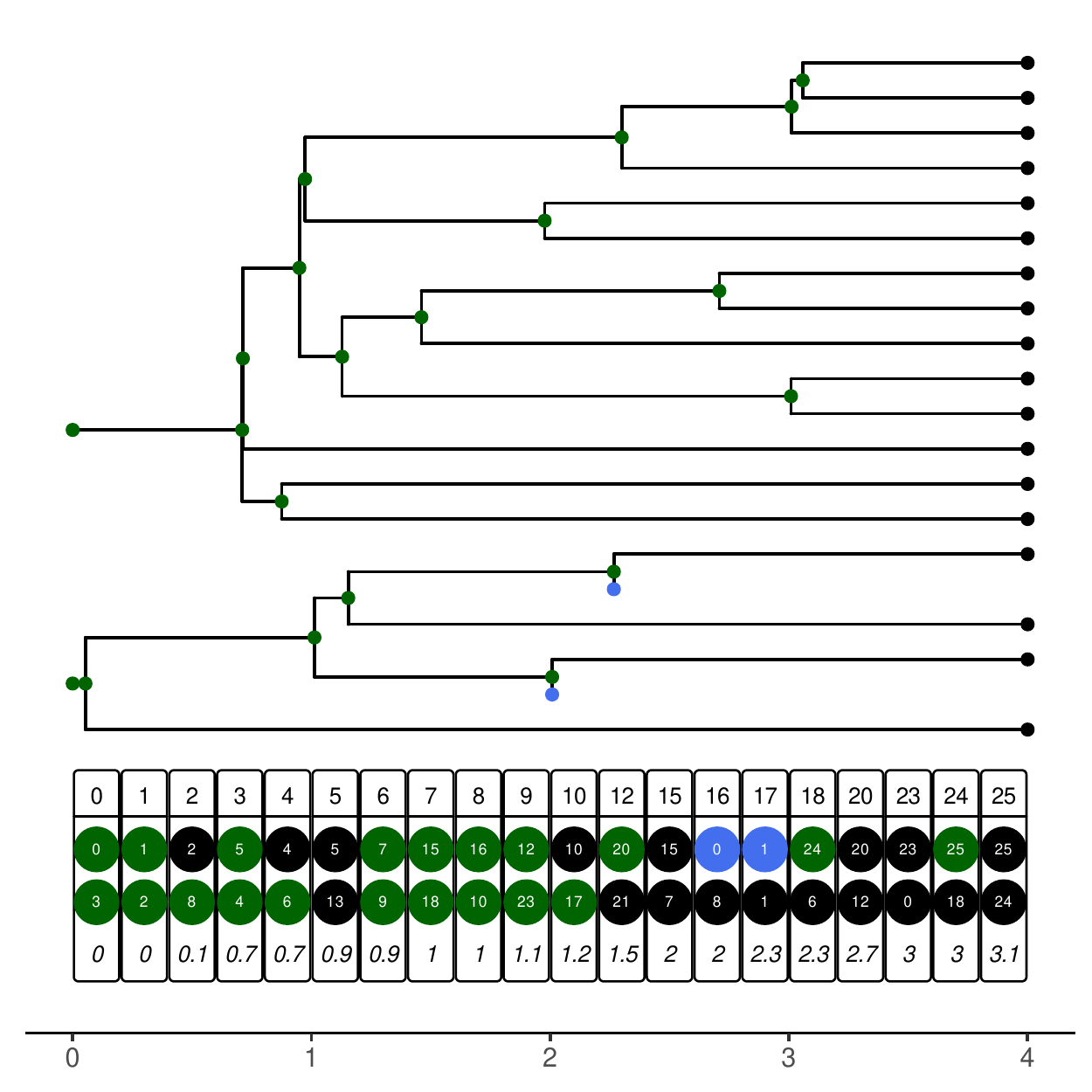} 

}

\end{knitrout}
  \caption{
    \label{fig:tree1}
    \textbf{Graphical and diagrammatic representations of a genealogy.}
    A genealogy can be represented graphically as a forest of 
    binary trees and diagrammatically as a sequence of nodes, as described in the text.
    In the node-sequence diagram, the nodes (depicted as narrow boxes) are ordered from left to right according to their time (the bottom number in each box).
    Each node has a name (actually a serial number, the top number in each box) and a \emph{pocket}, which holds two colored balls, one for each of the node's children.
    \emph{Leaves} are indicated with black points in the tree and black balls in the diagram;
    these represent living members of the population.
    \emph{Internal nodes} are indicated with green points in the tree and green balls in the diagram;
    these correspond to most recent common ancestors of subsets of the extant population.
    \emph{Samples} are indicated with blue points and balls.
    The horizontal axis is time.
    The genealogy depicted has two roots, i.e., two ancestors present at time 0, from whom all living members of the population are descended.
  }
\end{figure}

\paragraph{Genealogies}

A genealogy relates the members of a population that are alive at a given time.
Specifically, at any instant, every nonempty subset of living individuals has a most recent common ancestor.
Moreover, to each such subset is associated a unique time, viz., that at which the most distantly related lineages within the subset diverged.
Taken together, the collection of all such subsets and divergence times defines the genealogy, which is most commonly represented as a tree (\cref{fig:tree1}).
However, such representations are not unique and can be difficult to reason about.
We seek a representation that is unique and for which it will be straightforward to work out probabilistic properties.

Note that a genealogical tree has at least two kinds of nodes (\cref{fig:tree1}).
\emph{Leaves} represent members of the extant population at some time.
\emph{Internal nodes} represent ancestors, each of which is the most recent common ancestor of some subset of the extant population.
In addition, if the population is sampled, it is natural to represent the samples as nodes of third type.
The horizontal position of nodes in \cref{fig:tree1} is significant:
each node has an associated time.
In the case of leaves, the associated time is that of the extant population, i.e., that of the genealogy itself.

Hence, to characterize a genealogy, it is sufficient to note the time of the leaves and enumerate the internal nodes, noting for each its time and the identities of the nodes that descend immediately from it.
Since there are several kinds of nodes, we need some way of distinguishing between them.
To accomplish this, we introduce the notion of a \emph{colored ball}, which we define to be an ordered pair $(f,n)\in\Colors\times\N$, where $\Colors$ is a finite set.
We think of $f$ as the \emph{color} of the ball and $n$ as its \emph{name}.
Our convention is to associate black with leaves and green with internal nodes.
To handle sampling, we will need two additional colors, which we will take to be blue and red.
Thus $\Colors\coloneqq\{\green,\black,\blue,\red\}$.

We define a genealogical \emph{node} to be a triple $(n,t,w)$, where $n\in\N$ is the node's \emph{name}, $t\in\mathbb{R}$ is its \emph{time}, and $w$ is an (unordered) pair of colored balls, which we call the node's \emph{pocket}.
Given a node $\Node{p}$, we will use $n(\Node{p})$, $t(\Node{p})$, and $w(\Node{p})$ to denote the name, time, and pocket of $\Node{p}$, respectively.

A \emph{genealogy} is defined to be an ordered pair, $\Gen=\left(t,\left(\Node{p}_k\right)_{k=0}^{K-1}\right)$, where $t\in\Time$, $K\in\N$, the $\Node{p}_k$ are nodes, and the following conditions are satisfied:
\begin{enumerate}[(i)]
\item For all $j$, $t(\Node{p}_j)\le t(\Gen)$, i.e., no node time is later than that of the genealogy itself.
\item $j<k$ implies $t(\Node{p}_j)\le t(\Node{p}_k)$, i.e., the nodes are ordered in time.
\item $j\ne k$ implies $w(\Node{p}_j)\cap w(\Node{p}_k)=\emptyset$; the pockets of distinct nodes are disjoint.
\item $(\green,n(\Node{p}_j))\in w(\Node{p}_k)$ implies $j\ge k$; no green ball is held in the pocket of a later node.
\item For all $j$, $(\green,n(\Node{p}_j))\in w(\Gen)$; for every node, there is a green ball bearing the name of that node.
\end{enumerate}
Here, $w(\Gen)$ refers to the contents of the pockets in $\Gen$ collectively:
\begin{equation*}
  w(\Gen)\coloneqq\bigcup_{k=0}^{K-1} w(\Node{p}_k).
\end{equation*}
Following our usual convention, we use $t(\Gen)$, $K(\Gen)$, and $\Node{p}_k(\Gen)$ to refer to the time, the length, and the $k$-th node of genealogy $\Gen$, respectively.
We will use $\seq{\Gen}$ to refer to the node sequence of genealogy $\Gen$.
In a slight abuse of notation, we will write $\Node{p}\in\Gen$ when $\Node{p}$ is one node in $\seq{\Gen}$.

The black balls serve as pointers to members of the population extant at time $t$; the blue ones, to samples.
The green balls function as pointers to internal nodes.
In particular, a green ball held in the pocket of one node signifies that the node whose name matches that of the green ball is the immediate descendant of the first node.
Note that we allow a node to hold its own green ball, i.e., it is permissible that $(\green,n(\Node{p}))\in w(\Node{p})$.
Indeed, necessarily $(\green,n(\Node{p_0}))\in w(\Node{p_0})$ and, more generally, any node $\Node{p}$ with $(\green,n(\Node{p}))\in w(\Node{p})$ is a \emph{root} of the genealogy.
Note that nothing about our genealogy definition requires that the genealogical tree be connected:
multiple roots are allowed.

\begin{figure}

\begin{knitrout}\small
\definecolor{shadecolor}{rgb}{0.969, 0.969, 0.969}\color{fgcolor}

{\centering \includegraphics[width=0.7\linewidth]{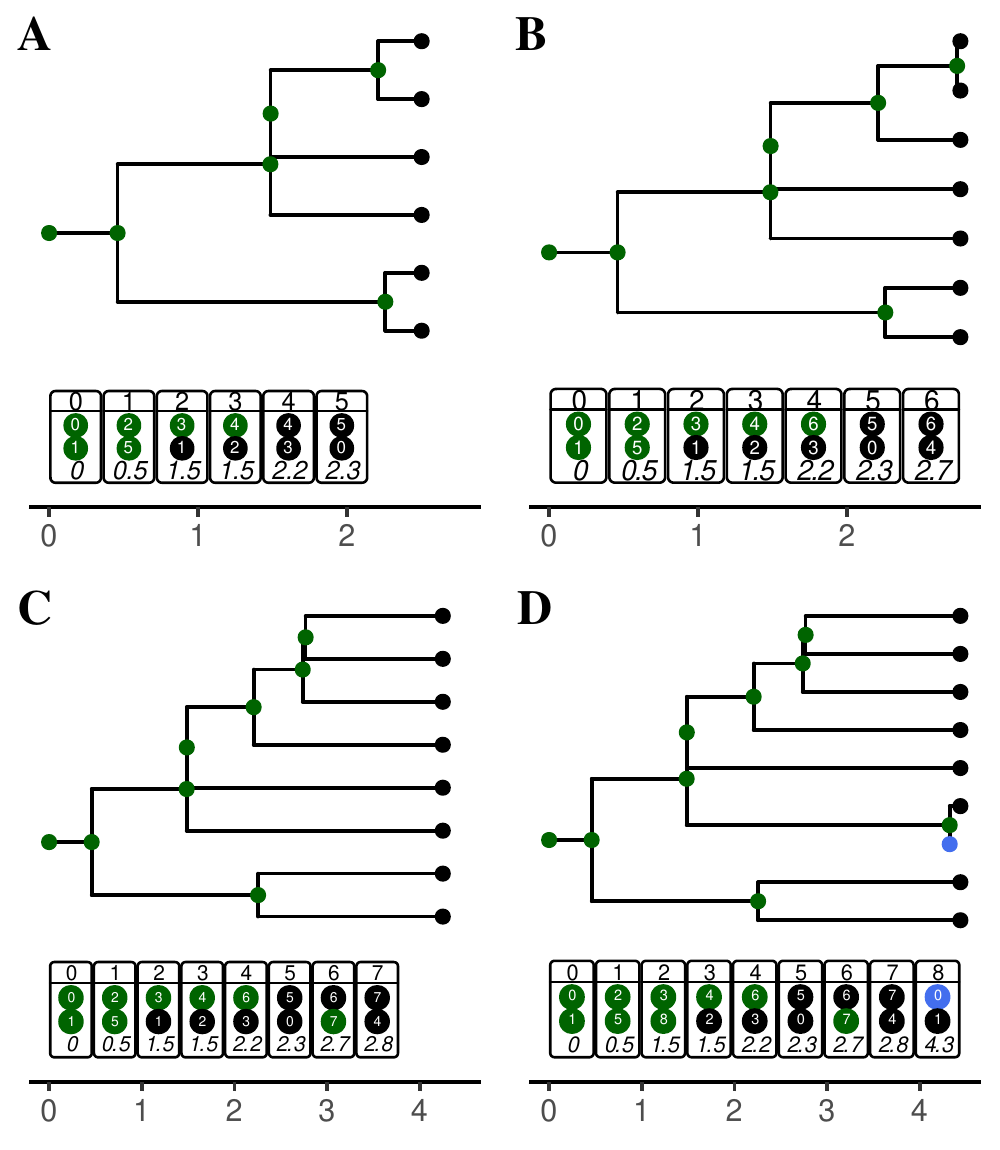} 

}

\end{knitrout}
  \caption{
    \label{fig:bds1}
    \textbf{Effects of births and sampling events on a genealogy.}
    Each panel shows the graphical representation of a genealogy as a tree, together with its representation as a node sequence.
    The horizontal axis is time.
    \textbf{(A--B)} The change from A to B shows the effect of a birth event on a genealogy.
    A new node (number 6) is introduced.
    The green ball corresponding to this node is exchanged for a randomly selected black ball (in this case, one held by node 4).
    \textbf{(C--D)} The change from C to D shows the effect of a sampling event on a genealogy.
    A new node (number 8), holding a blue ball, is introduced; its green ball is exchanged for a random black ball (in this case, the one held by node 2).
  }
\end{figure}

\begin{figure}
\begin{knitrout}\small
\definecolor{shadecolor}{rgb}{0.969, 0.969, 0.969}\color{fgcolor}

{\centering \includegraphics[width=0.7\linewidth]{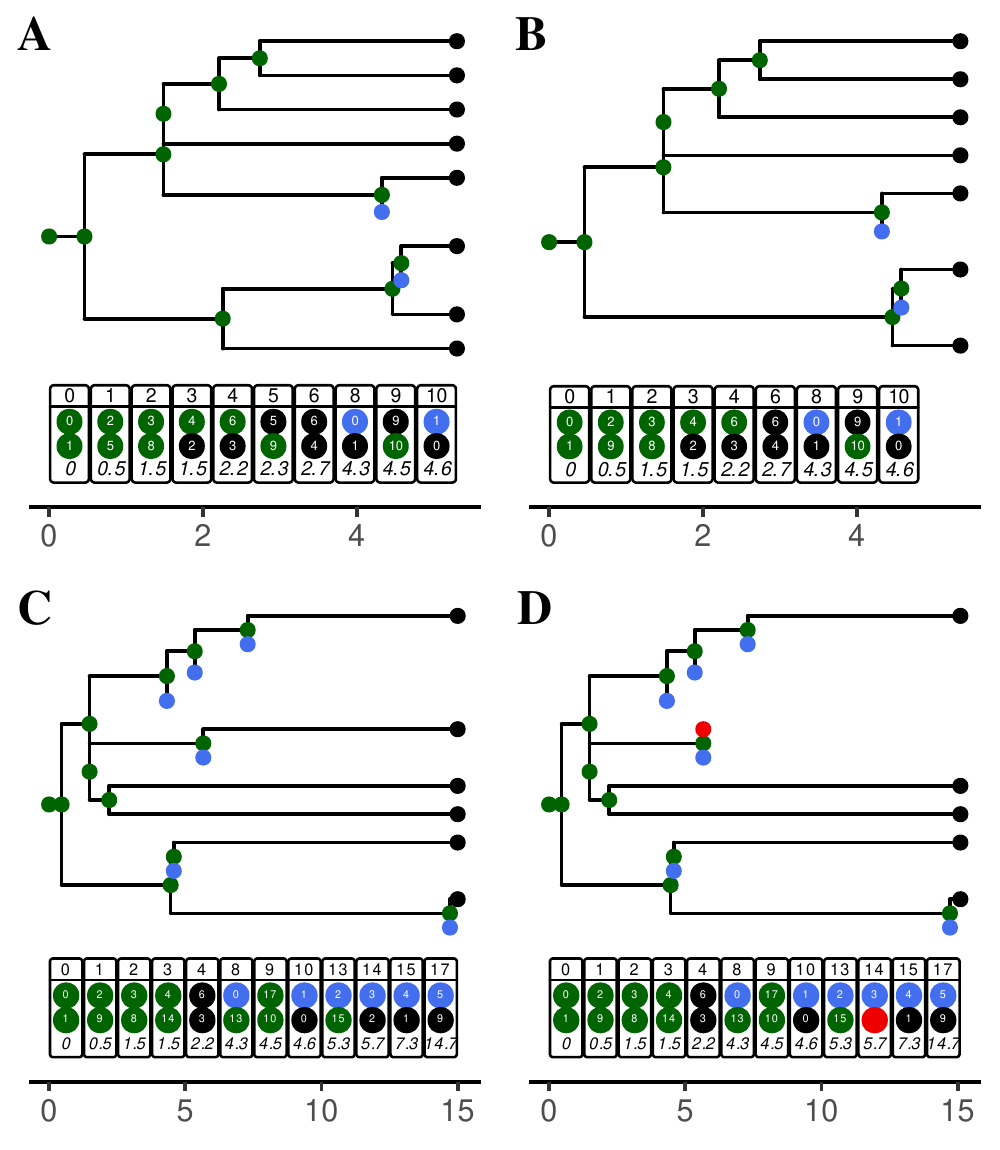} 

}

\end{knitrout}
  \caption{
    \label{fig:bds2}
    \textbf{Effect of deaths on a genealogy.}
    Each panel shows the graphical representation of a genealogy as a tree, together with its representation as a node sequence.
    The horizontal axis is time.
    \textbf{(A--B)} Passing from A to B, we see one possible effect of a death event on a genealogy.
    The random black ball chosen in this instance was the one held by node 5 in panel A.
    Since the other ball held by node 5 was not blue, node 5 is removed after first exchanging this other ball (green ball 9), for its own green ball (previously held by node 1).
    \textbf{(C--D)} The change from C to D shows the other possible effect of a death event.
    The random black ball chosen here was the one held by node 14.
    Since its other ball was blue, node 14 is not removed.
    Rather, it exchanges the selected black ball for a red ball.
    At any time, nodes with red balls represent samples on lineages that have since died out, while nodes holding one blue and one non-red ball represent samples on lineages that remain alive.
  }
\end{figure}

\paragraph{Effects of births, deaths, and sampling on a genealogy}

When births, deaths, or sampling events occur, these lead to changes in the genealogy (\cref{fig:bds1,fig:bds2}).
Here, we describe these changes in detail for each type of event in its turn.
In each instance, we assume $\Gen$ is a genealogy, as defined above, and that the event occurs at time $t$.
Each such event will involve one particular individual in the extant population.
Since the extant population is in 1-1 correspondence with the set of black balls in $w(\Gen)$, this amounts to choosing one black ball.
In the following, therefore, if $c_0<c_1<\dots<c_{m-1}$ are the names of the $m$ black balls in $w(\Gen)$ and $0\le n<m$, then $\Ball{b}=(\black,c_n)$ will be the selected black ball.

A birth in the population at time $t$ leads to the addition of a new internal node and a new leaf in the genealogy (\cref{fig:bds1}A--B).
Let $\Ball{b}$ be the $n$-th black ball, as above:
this is the parent of the newborn individual.
There is a unique node $\Node{p}\in\Gen$ such that $\Ball{b}\in w(\Node{p})=\{\Ball{b},\Ball{b}'\}$.
Let the name of the newborn be $c_m\coloneqq c_{m-1}+1$ and let $\Ball{g}=(\green,c_m)$ and $\Ball{b}''=(\black,c_m)$.
Construct a new node $\Node{p}'=(c_m,t,\{\Ball{g},\Ball{b}''\})$.
Now let node $\Node{p}'$ exchange $\Ball{g}$ with node $\Node{p}$ for $\Ball{b}$.
Thus, if before the swap we have $w(\Node{p})=\{\Ball{b},\Ball{b}'\}$ and $w(\Node{p}')=\{\Ball{g},\Ball{b}''\}$, after the swap we have $w(\Node{p})=\{\Ball{g},\Ball{b}'\}$ and $w(\Node{p}')=\{\Ball{b},\Ball{b}''\}$.
Finally, insert $\Node{p}'$ into the last position in the sequence of nodes.
Let $\mathrm{add}(\seq{\Gen},n)$ denote the resulting sequence of nodes.
\Cref{fig:bds1}A--B illustrates. 

A sample in the population at time $t$ results in the addition of one new internal node equipped with a new blue ball (\cref{fig:bds1}C--D).
Let $\Ball{b}$ be the $n$-th black ball selected as above:
this will be the individual sampled.
As before, there is a unique node $\Node{p}\in\Gen$ such that $\Ball{b}\in w(\Node{p})=\{\Ball{b},\Ball{b}'\}$.
Again, we construct a new node by taking $c_m=c_{m-1}+1$ and letting $\Node{p}'=(c_m,t,\{\Ball{g},\Ball{b}''\})$, where $\Ball{g}=(\green,c_m)$ and $\Ball{b}''=(\blue,q)$.
Again, we swap $\Ball{b}$ for $\Ball{g}$ between nodes $\Node{p}$ and $\Node{p}'$ and insert $\Node{p}'$ at the last position of the node-sequence.
Here we take the name, $q$, of the new blue ball to be the ordinal number of the sample, $q=|\{\Ball{b}\in w(\Gen)\;|\;\Ball{b}\ \text{is blue}\}|$. 
We denote the resulting sequence of nodes by $\mathrm{sample}(\seq{\Gen},n)$.
\Cref{fig:bds1}C--D illustrates. 

A death in the population at time $t$ can lead to the loss of one leaf and one internal node (\cref{fig:bds2}A--B).
The genealogy thus drops all record of the existence of the deceased individual.
On the other hand, samples represent recorded events: we do not wish to lose track of them.
Therefore, when a death would delete a sample, we prevent this from occurring using a red ball (\cref{fig:bds2}C--D).
To be precise, let $\Ball{b}$ be the $n$-th black ball selected as above: this is the individual who will die.
As usual, there is a unique node $\Node{p}\in\Gen$ such that $\Ball{b}\in w(\Node{p})=\{\Ball{b},\Ball{b}'\}$.
Let $\Ball{g}=(\green,n(\Node{p}))$ and let the unique node holding $\Ball{g}$ be denoted $\Node{p}'$.
If $\Ball{b}'$ is black, we swap $\Ball{b}'$ for $\Ball{g}$ and then delete $\Node{p}$ from the node sequence.
If $\Ball{b}'$ is blue, we replace $\Ball{b}$ with a red ball, with name matching that of $\Ball{b}'$, leaving everything else intact.
We use $\mathrm{drop}(\seq{\Gen},n)$ to denote the resulting sequence of nodes.
The two possible effects of sampling are illustrated in \Cref{fig:bds2}.

It is straightforward to verify that whenever $\Gen$ is a genealogy and $n<m$, $(t(\Gen),\mathrm{add}(\seq\Gen,n))$, $(t(\Gen),\mathrm{drop}(\seq\Gen,n))$, and $(t(\Gen),\mathrm{sample}(\seq\Gen,n))$ as just defined are all valid genealogies.
Note also that the $\mathrm{add}$ and $\mathrm{drop}$ procedures mirror their counterparts for the inventory process.
It follows that if $\Inv$ is an inventory, $\Gen$ a genealogy, and $\Inv$, $\Gen$ have the relation that $(\black,c)\in w(\Gen)$ if and only if $c\in\Inv$, then the same relation holds between $\mathrm{add}(\Inv)$ and $\mathrm{add}(\seq\Gen,n)$, $\mathrm{drop}(\Inv,n)$ and $\mathrm{drop}(\seq\Gen,n)$, and $\Inv$ and $\mathrm{sample}(\seq\Gen,n)$, respectively, for every $n<m$.

\paragraph{Genealogical event times}

Given a genealogy $\Gen$, the set of \emph{genealogical event times}, $\Eset(\Gen)$, is the set of all node times.
Several of its subsets are of interest.
In particular, we define
\begin{equation*}
  \begin{aligned}
    \Eset(\Gen) &\coloneqq \{t(\Node{p})\;|\;\Node{p}\in\Gen\}, \\
    \Masterset(\Gen) &\coloneqq \{t(\Node{p})\;|\;\Node{p}\in\Gen, w(\Node{p})\ \text{contains a green ball}\}, \\
    \Cset(\Gen) &\coloneqq \{t(\Node{p})\;|\;\Node{p}\in\Gen, w(\Node{p})\ \text{contains two green balls}\}, \\
    \Lset(\Gen) &\coloneqq \{t(\Node{p})\;|\;\Node{p}\in\Gen, w(\Node{p})\ \text{contains no green balls}\}, \\
    \Tset(\Gen) &\coloneqq \{t(\Node{p})\;|\;\Node{p}\in\Gen, w(\Node{p})\ \text{contains a blue ball}\},\\
    \Qset(\Gen) &\coloneqq \Masterset(\Gen) \cap \Tset(\Gen).
  \end{aligned}
\end{equation*}
With these definitions,
$\Masterset(\Gen)$ comprises the internal node times,
$\Cset(\Gen)$ is the set of branch times,
$\Tset(\Gen)$ holds the sample times,
and $\Qset(\Gen)$ is the set of \emph{direct-descent times},
i.e., the times of samples that are themselves directly ancestral to other samples.

\paragraph{Genealogy process}

We now proceed to define the \emph{genealogy process}, $\Gen_t$.
For $t\in\Time$ and $\omega\in\Omega$, we define $\Gen_t(\omega)=\left(t,\mathrm{geneal}(\restrict{\omega}{t})\right)$, where $\mathrm{geneal}$ is a deterministic procedure defined recursively for $\omega\in\mathring\Omega$ by
\fontsize{10pt}{12pt}\selectfont
\begin{equation}\label{eq:Gdef}
  \mathrm{geneal}(\omega)\coloneqq
  \begin{cases}
    \left(k,0,\{(\green,k),(\black,k)\}\right)_{k=0}^{I(\X_0(\omega))-1}, &\text{if}\ K(\omega)=0;\\
    \mathrm{add}(\mathrm{geneal}(\predec{\omega}),N_{K(\omega)}(\omega)), &\text{if}\ K(\omega)>0\ \text{and}\ U_{K(\omega)}\in\Bset;\\
    \mathrm{drop}(\mathrm{geneal}(\predec{\omega}),N_{K(\omega)}(\omega)), &\text{if}\ K(\omega)>0\ \text{and}\ U_{K(\omega)}\in\Dset;\\
    \mathrm{sample}(\mathrm{geneal}(\predec{\omega}),N_{K(\omega)}(\omega)), &\text{if}\ K(\omega)>0\ \text{and}\ U_{K(\omega)}\in\Sset;\\
    \mathrm{geneal}(\predec{\omega}), &\text{otherwise}.
  \end{cases}
\end{equation}
\normalsize
Thus, $\mathrm{geneal}(\omega)$ is a well defined node sequence for every jump sequence $\omega\in\mathring\Omega$.
Note that we initialize the genealogy process with a collection of root nodes.
The graphical representation of such a genealogy is a forest of single-leaf trees, each of which has zero branch length.
At each birth, death, or sample event, we modify the genealogy accordingly.
Note that the aforementioned parallelism between the $\mathrm{add}$, $\mathrm{drop}$, and $\mathrm{sample}$ procedures ensures that there is a deterministic map, $\mathrm{ext}$, such that $\Inv_t(\omega)=\mathrm{ext}(\Gen_t(\omega))$ for all $\omega\in\Omega$.

\begin{figure}
\begin{knitrout}\small
\definecolor{shadecolor}{rgb}{0.969, 0.969, 0.969}\color{fgcolor}

{\centering \includegraphics[width=1\linewidth]{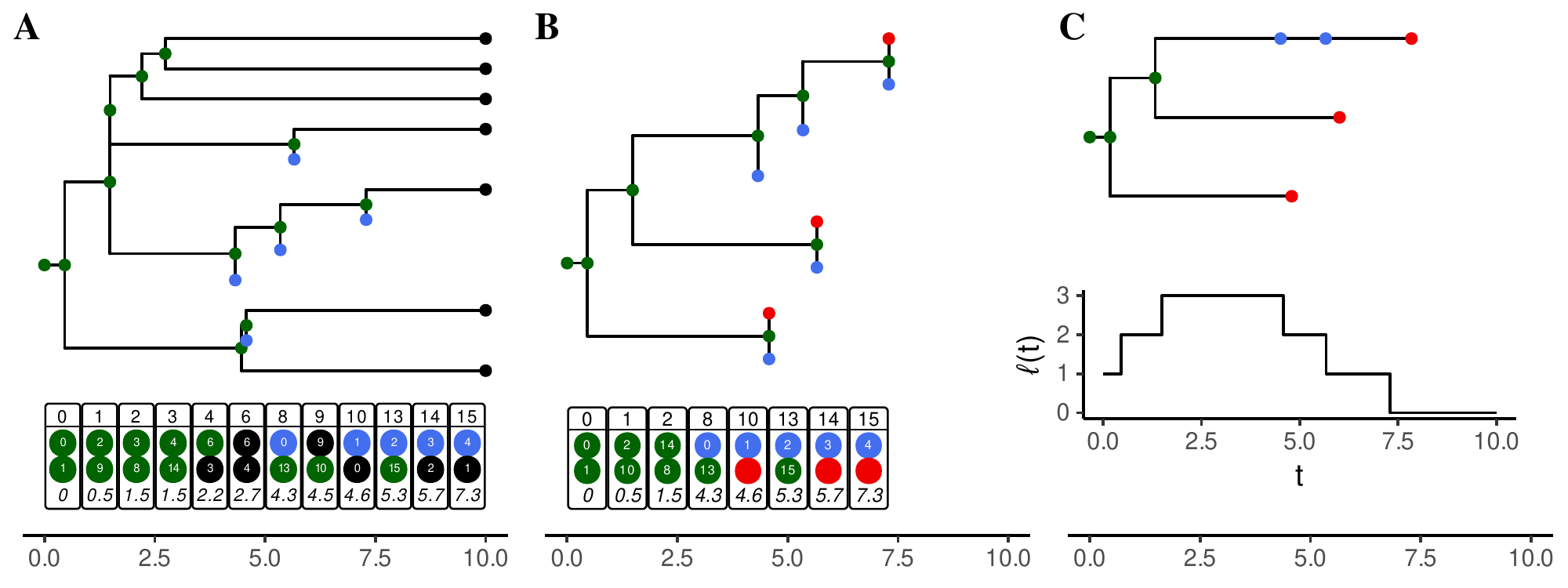} 

}

\end{knitrout}
  \caption{
    \label{fig:prune}
    \textbf{Pruning and the visible genealogy.}
    The genealogy depicted in \textbf{(A)} is pruned, i.e., all black balls are dropped according to the $\mathrm{drop}$ procedure described in the text.
    The resulting \emph{visible genealogy} is represented in \textbf{(B)}.
    A more compact graphical representation is displayed in \textbf{(C)}.
    In it, so-called \emph{red nodes}, i.e., those holding one red and one blue ball; are shown by red points in the graphical representation.
    Likewise \emph{blue nodes} (those holding one blue and one green ball) and \emph{green nodes} (those holding two green balls) are indicated by blue and green points, respectively.
    The inset shows the \emph{lineage count}, $\ell(t)$, the number of lineages present in the visible genealogy at time~$t$.
  }
\end{figure}

\paragraph{Visible genealogy process}

It is typically impossible to fully sample a population;
the genealogical relationships among unsampled lineages remain unobserved.
It is therefore of interest to study the genealogy that represents the relationships just among the samples.
This is most readily obtained from the full genealogy by a process of \emph{pruning}, which we proceed to describe.

Suppose $\Gen$ is a genealogy and let $\Inv=\mathrm{ext}(\Gen)$.
Let $\mathrm{obs}(\Gen)$ be the result of iteratively applying the $\mathrm{drop}$ operation defined above to $\Gen$ for each $c\in\Inv$.
Specifically, suppose $\Inv=\{c_0,\dots,c_{m-1}\}$.
Let $\Node{P}_0=\seq\Gen$ and $\Node{P}_{k}=\mathrm{drop}(\Node{P}_{k-1},c_{k-1})$, for $k=1,\dots,m$.
Then $\mathrm{obs}(\Gen)\coloneqq(t(\Gen),\Node{P}_{m})$.
\Cref{fig:prune} illustrates the pruning procedure.

For $\omega\in\Omega$, we define the \emph{visible genealogy process},
\begin{equation}
  \Vis_t(\omega)\coloneqq\textrm{obs}(\Gen_t(\omega)).
\end{equation}
Note that $\Vis_t$ so defined is a itself a genealogy.
\Cref{fig:prune}B depicts one example of $\Vis_t$ both graphically (as a tree) and diagrammatically (as a node-sequence).

\paragraph{Node color}

There are at most three kinds of nodes in a visible genealogy, distinguished by the contents of their pockets:
\begin{inparaenum}[(a)]
\item \emph{green nodes}, which have two green balls in their pockets;
\item \emph{blue nodes}, which have one blue and one green ball;
\item \emph{red nodes}, which have one red and one blue ball.
\end{inparaenum}
This distinction allows a compact graphical representation of the visible genealogy (\cref{fig:prune}C) and proves useful in deriving our main results as well.
Green nodes correspond to \emph{coalescence points} (branch points) in the visible genealogy;
red nodes, to leaves;
blue nodes correspond to direct-descent events, i.e., samples which are directly ancestral to other samples.
Thus, if $\Vis$ is a visible genealogy,
$\Eset(\Vis)=\Cset(\Vis)\cup\Qset(\Vis)\cup\Lset(\Vis)$, where
$\Cset(\Vis)$ is the set of times of the green nodes,
$\Qset(\Vis)$ contains the times of the blue nodes, and
$\Lset(\Vis)$ holds the times of the red nodes.
Moreover, $\Tset(\Vis)=\Qset(\Vis)\cup\Lset(\Vis)$ is the set of sample times and $\Masterset(\Vis)=\Qset(\Vis)\cup\Cset(\Vis)$ comprises the \emph{attachment times},
i.e., those points at which a sampled lineage attaches to the visible genealogy subtended by earlier samples.

\paragraph{Lineage count}

Given a genealogy, $\Vis$, at each time $t\in\Time$, there are a finite number, $\ell(t,\Vis)$, of lineages present in the genealogy at that time (\cref{fig:prune}C).
Evidently, one has
\begin{equation}\label{eq:elldef}
  \ell(t,\Vis)\coloneqq \sum_{e\in\Cset(\Vis)} \indicator{\halfopen{e,\infty}}(t) - \sum_{e\in\Lset(\Vis)} \indicator{\halfopen{e,\infty}}(t)=\sum_{e\in\Masterset(\Vis)} \indicator{\halfopen{e,\infty}}(t) - \sum_{e\in\Tset(\Vis)} \indicator{\halfopen{e,\infty}}(t).
\end{equation}
With this definition, note that $\ell$ is a \cadlag function of $t$.

\begin{figure}
\begin{knitrout}\small
\definecolor{shadecolor}{rgb}{0.969, 0.969, 0.969}\color{fgcolor}

{\centering \includegraphics[width=1\linewidth]{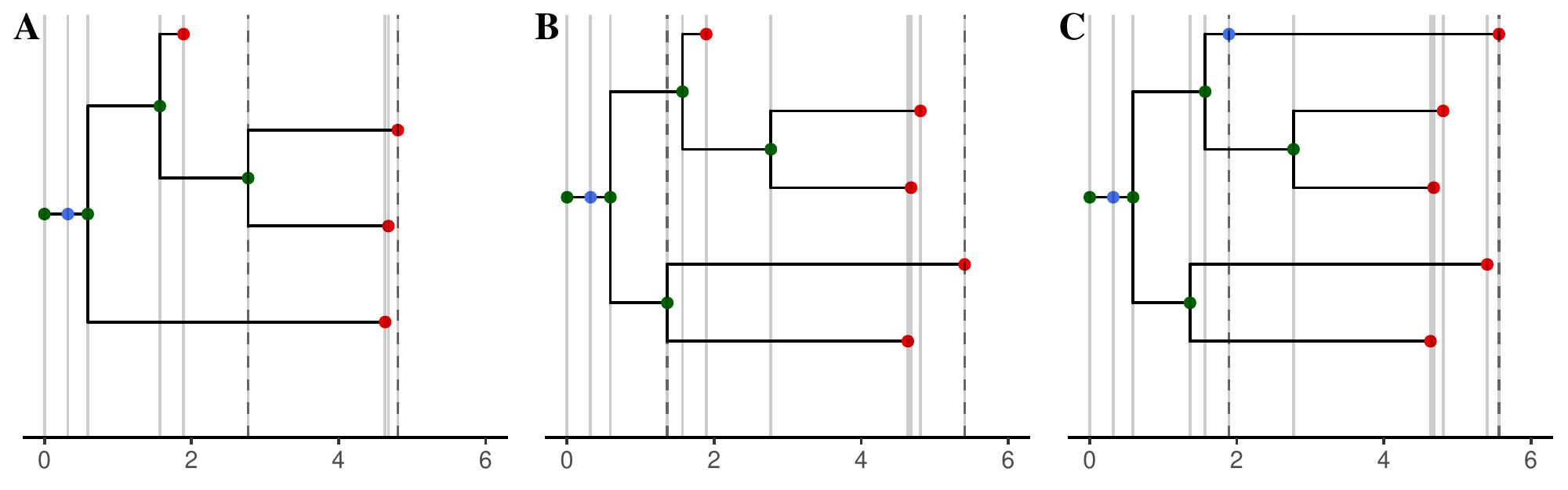} 

}

\end{knitrout}
  \caption{
    \label{fig:vischain}
    \textbf{Embedded chain of the visible genealogy process.}
    Three successive states of the embedded chain, $\Wch_i$, are shown.
    In panel A, the visible genealogy $\Wch_{5}$ is shown graphically.
    Note there are 5 samples represented (four red points and one blue one).
    The grey vertical lines indicate the genealogy event times, $\Eset(\Wch_5)$.
    Panels B and C show $\Wch_6$ and $\Wch_7$, respectively.
    In panel B, the new sample is depicted as the second red point from the bottom;
    it attaches to $\Wch_5$ at a green node (a coalescence point).
    In panel C, the new sample is the topmost red point; it attaches to $\Wch_6$ at a blue node (a direct-descent event).
    In each panel, the sample and attachment times of the latest sample are indicated by dashed vertical lines.
  }
\end{figure}

\paragraph{Embedded chain of the visible genealogy process}

Now, for $\omega\in\Omega$, consider the visible genealogy process $\Vis_t(\omega)$.
The sample times, $S_i$, form an increasing sequence.
Let $\Wch_i(\omega)\coloneqq\Vis_{S_i(\omega)}(\omega)$ be the embedded chain of the visible genealogy process.
Each genealogy in this chain builds on the previous one precisely in that one additional lineage is added;
\cref{fig:vischain} illustrates.
The terminal point of the new lineage is the latest sample time;
it attaches to the preceding genealogy at a random attachment time.
Let $A_i$ be the attachment time of the latest lineage in $\Wch_i$.
Note that the embedded chain is trivially Markov, since each $\Wch_i$ contains $\Wch_j$ within it, for $j<i$.

The proof of the following is immediate.
\begin{lemma}\label{lemma:ell}
  Let $s_i=\Tset(\Wch_i)\setminus\Tset(\Wch_{i-1})$ and $a_i=\Masterset(\Wch_i)\setminus\Masterset(\Wch_{i-1})$ be the sample and attachment times, respectively, of the $i$-th sample lineage.
  Then, if it is understood that $\ell(t,\Wch_0)=0$,
  \begin{equation*}
    \ell(t,\Wch_i)=\ell(t,\Wch_{i-1})+\indicator{\halfopen{a_i,s_i}}(t), \qquad \forall i>0.
  \end{equation*}
\end{lemma}

\section{Results}
\label{sec:results}

Let $P_{\Wch_i|\Hist}$ denote the probability density of $\Wch_i$ conditional on $\Hist_{S_i}$.
The Markovity of $\Wch_i$ gives us
\begin{equation*}
  P_{\Wch_i|\Hist}=P_{\Wch_1|\Hist}\,\prod_{j=2}^iP_{\Wch_j|\Wch_{j-1},\Hist}.
\end{equation*}
Now, conditional on $\Hist_{S_1}$, $\Wch_1$ consists a.s.\ of one red node and one root node.
Thus $P_{\Wch_1|\Hist}=1$.

Let us compute $P_{\Wch_j|\Wch_{j-1},\Hist}$.
We fix $A_j=a_j$, $\Wch_j=w_j$, $S_j=s_j$, and $\Hist_{s_i}=h=\left(s_i,\left(t_k,u_k\right)_{k=0}^K\right)$.
We write $x_k=\X_{t_k}$, and define the sets
$\Omega_j\coloneqq\{\omega\in\Omega\;|\;\Wch_{j}(\omega)=w_j,\Hist_{s_i}(\omega)=h\}\subset\Omega$ and
$\Omega_{jk}\coloneqq N_k(\Omega_j)\subset\N$, for $j\le i$ and $k\le K$.
There is a bijection between $\Omega_j$ and $\prod_{k=1}^K\Omega_{jk}$.
The conditional independence of the $N_k$ (\cref{eq:condH}) allows us to write
\begin{equation*}
  \begin{aligned}
    P_{\Wch_j|\Wch_{j-1},\Hist}(w_j|w_{j-1},h)&=\sum_{\omega\in\Omega_j}P_{\Master|\Hist}(\Master(\omega)|h)
    =\sum_{\omega\in\Omega_j}\prod_{k=1}^K\auxmeas_{u_k,x_k-u_k}(N_k(\omega))\\
    &=\prod_{k=1}^K\sum_{n_k\in\Omega_{jk}}\auxmeas_{u_k,x_k-u_k}(n_k).
  \end{aligned}
\end{equation*}
Since the $\auxmeas$ are uniform, it remains only to count up the number of elements in the sets $\Omega_{jk}$, i.e., the number of choices of $N_{k}$ that are consistent with $\Wch_j$, for each $j$ and $k$.

Define $q_{jk}\coloneqq\sum_{n_k\in\Omega_{jk}}\auxmeas_{u_k,x_k-u_k}(n_k)$ and,
for the moment, let $I_k=I(x_k)$, $\ell_{jk}=\ell(t_k,\Wch_{j-1})$.
Note that $\ell_{1k}=0$ for all $k$.
There are several cases to consider:
\begin{enumerate}[(a)]
\item If $t_k\notin\halfopen{a_j,s_j}$, or if $u_k\notin\Bset\cup\Sset$, then all choices $n_k$ are compatible, so we have $q_{jk}=1$.
\item If $t_k\in\Masterset(\Wch_{j-1})$, then again, all choices $n_k$ are compatible and $q_{jk}=1$.
\item If $t_k\in(a_j,s_j)\setminus\Masterset(\Wch_{j-1})$ and $u_k\in\Sset$, then there was a sample event at time $t_k$ but the $j$-th sample lineage did not directly descend from it.
  Since $\Sset\cap(\Bset\cup\Dset)=\emptyset$, $\Inv(t_k-)=\Inv(t_k)$, i.e., the inventory of the population did not change at time $t_k$.
  There were $I_k$ individuals at this time.
  However, $\ell_{jk}$ of these could not have been chosen for the sample (or we would have $t_k\in\Masterset(\Wch_{j-1})$).
  Only one of these was of the $j$-th sample lineage.
  Therefore $q_{jk}=1-{1}/({I_k-\ell_{jk}})$ in this case.
\item If $t_k=a_j$ and $u_k\in\Sset$, then the $j$-th lineage attaches to $\Wch_{j-1}$ in a direct-descent event.
  There were $I_k-\ell_{jk}$ individuals who might have been sampled, just one of whom was of the $j$-th lineage.
  Therefore, $q_{jk}=1/(I_k-\ell_{jk})$ in this case.
\item \label{item:bNoC} If $t_k\in(a_j,s_j)\setminus\Masterset(\Wch_{j-1})$ and $u_k\in\Bset$, then there was a birth event at time $t_k$ but the $j$-th sample lineage was not involved in it.
  Thus the population size just before $t_k$ was $I_k-1$.
  The node newly added has a pair of black balls, one corresponding to the newborn, the other to the parent; these are indistinguishable.
  There are $\binom{I_k}{2}$ such pairs, but $\binom{\ell_{jk}}{2}$ of these can be excluded because none of the $\ell_{jk}$ lineages of $\Wch_{j-1}$ coalesces here.
  Of the remaining pairs, exactly $\ell_{jk}$ involve the sole individual of the $j$-th sample lineage and one individual of one of the lineages in $\Wch_{j-1}$.
  Therefore, $q_{jk}=1-\ell_{jk}/(\binom{I_k}{2}-\binom{\ell_{jk}}{2})$ in this case.
\item Finally, we treat the case $t_k=a_j$, $u_k\in\Bset$.
  The logic of case \ref{item:bNoC} applies, but here the coalescence event did occur.
  Precisely one of the $\binom{I_k}{2}-\binom{\ell_{jk}}{2}$ pairs that might have coalesced is the one consisting of the individual of the $j$-th sample lineage and the individual in the sample lineage to which the $j$-th sample lineage did attach.
  Therefore, $q_{jk}=1/(\binom{I_k}{2}-\binom{\ell_{jk}}{2})$ in this case.
\end{enumerate}

To summarize, we have
\begin{equation}\label{eq:qcases}
  q_{jk}=
  \begin{cases}
    1, &\text{if}\ t_k\notin\halfopen{a_j,s_j},\\
    1, &\text{if}\ u_k\notin\Bset\cup\Sset,\\
    1, &\text{if}\ t_k\in\Cset(\Wch_{j-1})\cup\Qset(\Wch_{j-1}),\\[12pt]
    1-\frac{1}{I_k-\ell_{jk}}, &\text{if}\ t_k\in(a_j,s_j)\setminus\Qset(\Wch_{j-1})\ \text{and}\ u_k\in\Sset,\\[12pt]
    \frac{1}{I_k-\ell_{jk}}, &\text{if}\ t_k=a_j\in\Lset(\Wch_{j-1}),\\[12pt]
    1-\frac{\ell_{jk}}{\binom{I_k}{2}-\binom{\ell_{jk}}{2}}, &\text{if}\ t_k\in(a_j,s_j)\ \text{and}\ u_k\in\Bset,\\[12pt]
    \frac{1}{\binom{I_k}{2}-\binom{\ell_{jk}}{2}}, &\text{if}\ t_k=a_j\in\Cset(\Wch_j).
  \end{cases}
\end{equation}

This establishes the first of our main results, namely:

\begin{thm}\label{thm:hardway}
  With the definitions as above, $P_{\Wch_i|\Hist}(w|h) = \prod_{j=1}^{i}\prod_{k=1}^{K}q_{jk}$.
\end{thm}

Although this result is compactly stated, the computation it suggests is awkward.
For example, to compute $P_{\Wch_i|\Hist}(w|h)$ via Monte Carlo, one must simulate and store the entire history $\Hist_t$, in order to compute each of the $\prod_{k}q_{jk}$.
For all but quite small populations, this will prove impracticable.
Simply exchanging the order of the product in \cref{thm:hardway}, however, yields a scheme which requires only simulation and temporary storage of the population process $\X_t$.

\begin{thm}\label{thm:easyway}
  Suppose $\Vis_t$ is the visible genealogy process induced by $\X_t$ and $\Hist_t$ is the history process.
  Fix $\Hist_t=h=\left(t,\left(t_k,u_k\right)_{k=0}^{K}\right)$ and let $x_t=\sum_{k=0}^Ku_k\,\indicator{\halfopen{t_k,\infty}}(t)$.
  Let $\Uset(h)$ be the set of unobserved births in history $h$, i.e., $\Uset(h)=\{t_k(h)\;|\;k>0\ \text{and}\ u_k(h)\in\Bset\}\setminus\Cset(\Vis_t)$.
  Then
  \begin{equation*}
    P_{\Vis_t|\Hist_t}(\Vis_t|h) = \frac{\prod_{e\in \Uset(h)}\left(1-\frac{\binom{\ell(e,\Vis_t)}{2}}{\binom{I(x_e)}{2}}\right)\,\prod_{e\in\Lset(\Vis_t)}\left(1-\frac{\ell(e,\Vis_t)}{I(x_e)}\right)}{\prod_{e\in\Cset(\Vis_t)}\binom{I(x_e)}{2}\,\prod_{e\in\Qset(\Vis_t)}I(x_e)}.
  \end{equation*}
\end{thm}
\begin{proof}
  $\Vis_t=\Wch_i$ for some $i$.
  We compute $Q_k\coloneqq\prod_{j=1}^iq_{jk}$ for each $k$.
  As above, let $I_k=I(x_k)$, $\ell_{jk}=\ell(t_k,\Wch_{j-1})$.

  For the first three cases of \cref{eq:qcases}, we have $Q_k=1$.

  Suppose $t_k\in\Lset(\Wch_i)$.
  Then $u_k\in\Sset$, and we have
  \begin{equation*}
    Q_k=\prod_{j=1}^{i}\frac{I_k-\ell_{jk}-1}{I_k-\ell_{jk}}=\frac{I_k-\ell_{ik}-1}{I_k-\ell_{1k}}\,\cdot\,\prod_{j=1}^{i-1}\frac{I_k-\ell_{jk}-1}{I_k-\ell_{j+1,k}}.
  \end{equation*}
  By \cref{lemma:ell}, $\ell_{j+1,k}=\ell_{jk}+1$, which implies
  \begin{equation*}
    Q_k=\frac{I_k-\ell_{ik}-1}{I_k-\ell_{1k}}=1-\frac{\ell_{i+1,k}}{I_k}.
  \end{equation*}

  Now consider $t_k\in\Qset(\Wch_i)$.
  Let $m$ be the unique integer such that $t_k\in\Qset(\Wch_m)\setminus\Qset(\Wch_{m-1})$.
  Then we have
  \begin{equation*}
    q_{jk}=
    \begin{cases}
      1-\frac{1}{I_k-\ell_{jk}}, &j<m,\\[6pt]
      \frac{1}{I_k-\ell_{mk}}, &j=m,\\[6pt]
      1, &j>m.
    \end{cases}
  \end{equation*}
  It follows that
  \begin{equation*}
    Q_k=\frac{1}{I_k-\ell_{mk}}\,\cdot\,\prod_{j=1}^{m-1}\frac{I_k-\ell_{jk}-1}{I_k-\ell_{jk}}=\frac{1}{I_k}.
  \end{equation*}

  Now consider the case $t_k\notin\Cset(\Wch_i)$ but $u_k\in\Bset$.
  Again using \cref{lemma:ell}, we have
  \begin{equation*}
    Q_k=\prod_{j=1}^{i}\left(1-\frac{\ell_{jk}}{\binom{I_k}{2}-\binom{\ell_{jk}}{2}}\right)=1-\frac{\binom{\ell_{i+1,k}}{2}}{\binom{I_k}{2}}.
  \end{equation*}

  Finally, suppose $t_k\in\Cset(\Wch_i)$.
  Let $m$ be the unique integer such that $t_k\in\Cset(\Wch_m)\setminus\Cset(\Wch_{m-1})$.
  Then
  \begin{equation*}
    q_{jk}=
    \begin{cases}
      1-\frac{\ell_{jk}}{\binom{I_k}{2}-\binom{\ell_{jk}}{2}}, &j<m,\\[10pt]
      \frac{1}{\binom{I_k}{2}-\binom{\ell_{jk}}{2}}, &j=m\\[10pt]
      1, &j>m.
    \end{cases}
  \end{equation*}
  Once again, a straightforward calculation yields
  \begin{equation*}
    Q_k=\frac{1}{\binom{I_k}{2}}.
  \end{equation*}
  
\end{proof}

\paragraph{The unnormalized nonlinear filter (DMZ) equation}

\Cref{thm:easyway} gives us an expression for the likelihood of $\Vis_t$ given the history $\Hist_t$.
We now seek to integrate out the dependence on $\Hist_t$.
The form of \cref{thm:easyway} implies that this can be done in a sequential fashion, working from earlier times to later ones, thus avoiding the need either to work backward in time or to store the full history process.
Indeed, to compute $P_{\Vis_t|\Hist_t}$, we progressively accumulate factors, one for each event in the history and one for each event in the visible genealogy.
Each such term depends only on the state of the population process, $\X$ at the time of that event.
We can therefore integrate out the history by integrating over the possible values of $\X$ at each time.

Given a visible genealogy $\Vis$, and a time $0\le t\le t(\Vis)$, let us define the \emph{partial genealogy}, $\restrict{\Vis}{t}$, to be that portion of the visible genealogy lying to the left of $t$.
For each $t<t(\Vis)$, we will define the \emph{partial weight}, $w(t,x,\Vis)$, in such a way that
\begin{equation}\label[pluralequation]{eq:wgoal}
  \begin{gathered}
    w(t,x,\Vis) = P_{\Vis_t|\X_t}(\Vis\:|\:x) \qquad \text{and} \qquad
    \sum_{x\in\Xspace} w(t,x,\Vis) = P_{\Vis_t}(\Vis).
  \end{gathered}
\end{equation}
Note that for $t<t(\Vis)$, $w(t,x,\Vis)$ is not itself a likelihood of any subset of the data.
However, when $t=t(\Vis)$, the sum (over $x$) of the partial weights will equal the likelihood of the genealogy, unconditional on the history process.
The partial weights are defined by an initial-value problem that they must satisfy.
We proceed to derive this now.

Suppose $\Vis$ is a visible genealogy, $0\le t< t(\Vis)$, and $\delta{t}>0$.
Within the interval $\halfopen{t,t+\delta{t}}$, the following events are exhaustive and mutually exclusive:
\begin{inparaenum}[(a)]
\item\label[inpar]{it:nada} nothing occurred,
\item\label[inpar]{it:more} more than one event occurred,
\item\label[inpar]{it:other} an event occurred which was neither a birth nor a sample, and no other event occurred,
\item\label[inpar]{it:nocoal} a birth event, and no other, occurred, and this birth event was not a coalescence event in $\Vis$,
\item\label[inpar]{it:coal} a birth event, which was also a coalescence event, occurred, and no other,
\item\label[inpar]{it:dead} a sample event, and no other, occurred, and this sample was a direct-descent event in $\Vis$,
\item\label[inpar]{it:live} a sample event, which was not a direct-descent event, occurred, and no other.
\end{inparaenum}
Now, if $t\notin\Eset(\Vis)$, our non-explosion assumption implies that we can choose $\delta{t}$ sufficiently small so that $\Eset(\Vis)\cap\halfopen{t,t+\delta{t}}=\emptyset$.
In this case, the only possible events are \cref{it:nada,it:more,it:other,it:nocoal}.
Accordingly, we desire that
\begin{equation}\label{eq:wdiff1}
  \begin{aligned}
    w(t+\delta{t},x,\Vis)=
    &\left(1-\sum_{u\in\Uspace}\jumprate_u(t,x)\,\delta{t}\right)\,w(t,x,\Vis)\\
    &+\sum_{\mathclap{u\in\Uspace\setminus\Bset\setminus\Sset}}\jumprate_u(t,x-u)\,w(t,x-u,\Vis)\,\delta{t}\\
    &+\sum_{u\in\Bset}\jumprate_u(t,x-u)\,\left(1-\frac{\binom{\ell(t,\Vis)}{2}}{\binom{I(x)}{2}}\right)\,w(t,x-u,\Vis)\,\delta{t}
    +o(\delta{t}).
  \end{aligned}
\end{equation}
Rearranging \cref{eq:wdiff1} and taking $\delta{t}\downarrow 0$ in the usual way, we obtain, for $t\notin\Eset(\Vis)=\Cset(\Vis)\cup\Qset(\Vis)\cup\Lset(\Vis)$,
\begin{equation}\label{eq:wde1}
  \begin{aligned}
    \pd{w}{t}(t,x,\Vis)=&\sum_{\mathclap{u\in\Uspace}}\left[\jumprate_u(t,x-u)\,w(t,x-u,\Vis)-\jumprate_u(t,x)\,w(t,x,\Vis)\right]\\
    &-\sum_{u\in\Sset}\jumprate_u(t,x-u)\,w(t,x-u,\Vis)\\
    &-\sum_{u\in\Bset}\jumprate_u(t,x-u)\,\frac{\binom{\ell(t,\Vis)}{2}}{\binom{I(x)}{2}}\,w(t,x-u,\Vis).\\
  \end{aligned}
\end{equation}
On the other hand, if $t\in\Eset(\Vis)=\Cset(\Vis)\cup\Qset(\Vis)\cup\Lset(\Vis)$, we desire that
\begin{equation}\label{eq:wdiff2}
  w(t,x,\Vis)=
  \begin{cases}
    \displaystyle\sum_{u\in\Bset}\frac{\jumprate_u(t,x-u)}{\mu}\,\frac{1}{\binom{I(x)}{2}}\,w(\prev{t},x-u,\Vis), &t\in\Cset(\Vis),\\[13pt]
    \displaystyle\sum_{u\in\Sset}\frac{\jumprate_u(t,x-u)}{\mu}\,\frac{1}{I(x)}\,w(\prev{t},x-u,\Vis), &t\in\Qset(\Vis),\\[13pt]
    \displaystyle\sum_{u\in\Sset}\frac{\jumprate_u(t,x-u)}{\mu}\,\left(1-\frac{\ell(t,\Vis)}{I(x)}\right)\,w(\prev{t},x-u,\Vis), &t\in\Lset(\Vis).\\
  \end{cases}
\end{equation}
Here, $\prev{t}$ indicates the left limit.

Making use of the Dirac delta function, $\delta(t)$, we can combine \cref{eq:wde1,eq:wdiff2} into a single equation, the analogue of the Duncan-Mortensen-Zakai equation \citep{Zakai1969} for this problem:
\begin{equation}\label{eq:dmz}
  \begin{aligned}
    \pd{w}{t}(t,x,\Vis)=&\sum_{u\in\Uspace}\left[\jumprate_u(t,x-u)\,w(t,x-u,\Vis)-\jumprate_u(t,x)\,w(t,x,\Vis)\right]\\
    &-\sum_{u\in\Sset}\jumprate_u(t,x-u)\,w(t,x-u,\Vis)\\
    &-\sum_{u\in\Bset}\jumprate_u(t,x-u)\,\frac{\binom{\ell(t,\Vis)}{2}}{\binom{I(x)}{2}}\,w(t,x-u,\Vis)\\
    &+\sum_{e\in\Cset(\Vis)}\delta(t-e)\,\left\{\sum_{u\in\Bset}\frac{\jumprate_u(t,x-u)}{\mu}\,\frac{1}{\binom{I(x)}{2}}\,w(t,x-u,\Vis)\right\}\\
    &+\sum_{e\in\Qset(\Vis)}\delta(t-e)\,\left\{\sum_{u\in\Sset}\frac{\jumprate_u(t,x-u)}{\mu}\,\frac{1}{I(x)}\,w(t,x-u,\Vis)\right\}\\
    &+\sum_{e\in\Lset(\Vis)}\delta(t-e)\,\left\{\sum_{u\in\Sset}\frac{\jumprate_u(t,x-u)}{\mu}\,\left(1-\frac{\ell(t,\Vis)}{I(x)}\right)\,w(t,x-u,\Vis)\right\}\\
    &-\sum_{e\in\Eset(\Vis)}\delta(t-e)\,w(t,x,\Vis).\\
  \end{aligned}
\end{equation}
The appearance in \cref{eq:dmz} of the rate, $\mu$, of the Poisson point process, the probability measure of which is the base measure for our probability densities, serves as a reminder that the numerical values of these densities depend on the choice of time unit.

For \cref{eq:dmz} to be valid, we must insist that $I(x)\ge\ell(t,\Vis)$.
Moreover, if $I(x)<\ell(t,\Vis)$, then the visible genealogy $\Vis$ is incompatible with $\X_t=x$.
Therefore, we put $w(t,x,\Vis)=0$, for all $x$ such that $I(x)<\ell(t,\Vis)$.
The proper initial condition is clearly
\begin{equation}\label{eq:wic}
  w(0,x,\Vis_t)=p_0(x).
\end{equation}
With these definitions, it is a straightforward matter to verify that the unique $w$ satisfying \cref{eq:wic,eq:dmz} also satisfies \cref{eq:wgoal}.
In particular, the likelihood of a visible genealogy $\Vis$ is
\begin{equation}\label{eq:lik}
  \lik_{\Vis}=\sum_x w(t(\Vis),x,\Vis).
\end{equation}

\section{Illustrative examples}
\label{sec:illus}

In this section, we return to some of the examples of \cref{sec:examples}.
For each one, we specialize \cref{eq:dmz} and perform likelihood calculations on simulated data.
Codes for the following (and for all the figures in the paper) are available as a \href{https://doi.org/10.5281/zenodo.5758900}{Zenodo digital archive}.

\paragraph{Linear birth-death-sampling process}

For any given model, \cref{eq:dmz,eq:wde1,eq:wdiff2} take specific forms.
In the case of the linear birth-death-sampling process (\cref{sec:examples}), if we write $w(t,n)$ for the partial weight associated with population size $n$, \Cref{eq:wde1} becomes
\begin{equation}\label{eq:wde1_lbdp}
  \frac{\partial{w}}{\partial{t}} = \lambda\,(n-1)\,\left(1-\frac{\binom{\ell}{2}}{\binom{n}{2}}\right)\,w(t,n-1)+\delta\,(n+1)\,w(t,n+1)-(\lambda+\delta+\psi)\,n\,w(t,n),
\end{equation}
which holds for $t\notin\Eset(\Vis)$.
We can integrate \cref{eq:wde1_lbdp} forward in time from each genealogical event to the next.
We then adjust $w$ according to the nature of the event, as follows:
\begin{equation}\label{eq:wdiff2_lbdp}
  w(t,n)=
  \begin{cases}
    \frac{\lambda\,(n-1)}{\mu}\,\frac{1}{\binom{n}{2}}\,w(\prev{t},n-1), &t\in\Cset(\Vis),\\[13pt]
    \frac{\psi}{\mu}\,w(\prev{t},n), &t\in\Qset(\Vis),\\[13pt]
    \frac{\psi}{\mu}\,\left(n-\ell\right)\,w(\prev{t},n), &t\in\Lset(\Vis).\\
  \end{cases}
\end{equation}
In these equations, it is understood that $w(t,n)=0$ for $n<\ell(t)$.

The astute reader will have noticed that in \cref{eq:wde1_lbdp,eq:wdiff2_lbdp}, the only state variable upon which $w$ depends is the population size $n$, while in \cref{sec:examples}, we specify a two-dimensional state space for the linear birth-death-sampling process with coordinates $(n,g)$, $g$ being the cumulative number of samples to time $t$.
Since, as is easily verified, solving \cref{eq:dmz} results in $w(t,n,g,\Vis)\ne 0$ if and only if $g$ is precisely equal to the number of samples in $\restrict{\Vis}{t}$, there is no need to keep track of the the dependence of $w$ on $g$.
The same simplification is used in the other examples we consider below.

\begin{figure}
\begin{knitrout}\small
\definecolor{shadecolor}{rgb}{0.969, 0.969, 0.969}\color{fgcolor}

{\centering \includegraphics[width=1\linewidth]{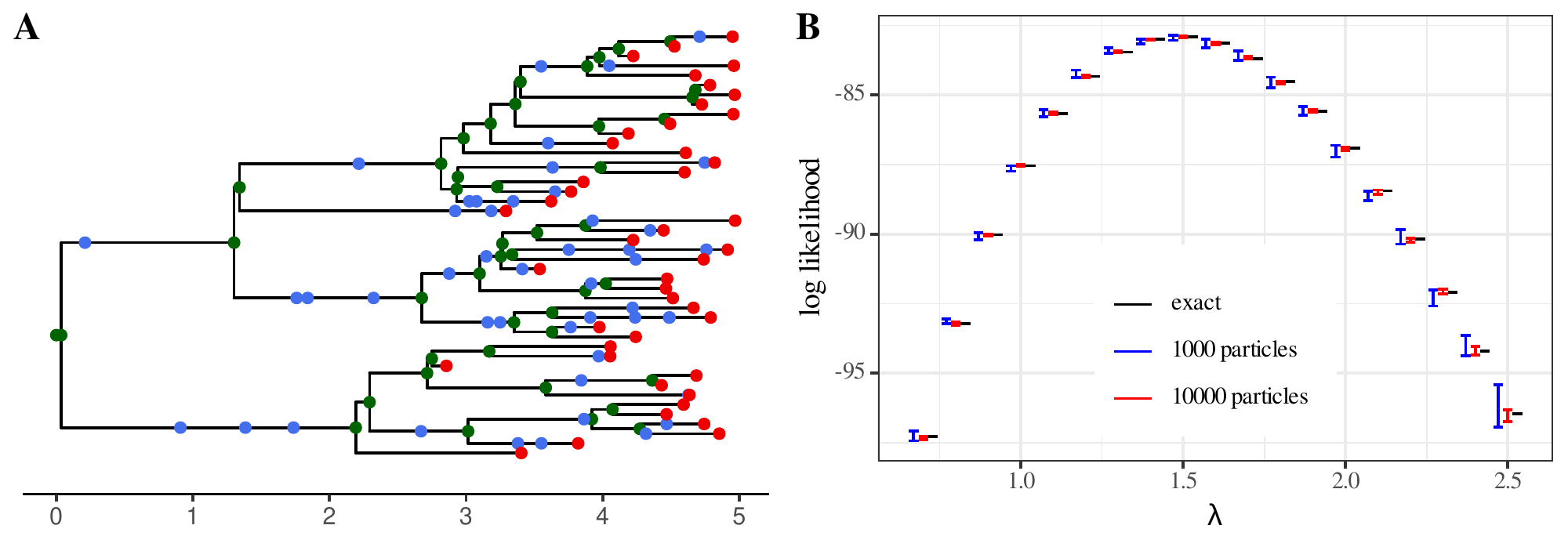} 

}

\end{knitrout}
  \caption{
    \textbf{Computing the likelihood for genealogies induced by the linear birth-death-sampling process.}
    \textbf{(A)} A simulated genealogy for $\lambda=1.5$, $\mu=0.8$, $\psi=1$.
    As usual, blue and red points correspond to samples;
    green points represent branch-points.
    \textbf{(B)} Change in the log likelihood as a function of position along a line passing through the true parameter in the $\lambda$-direction.
    The red and blue error bars show the estimates obtained using the particle filter (mean ${\pm}2$~s.e.), with different amounts of effort (i.e., number of particles);
    the black shows the exact log likelihood, which is available in closed form in this case.
    With increasing computational effort, the estimates converge on the exact value.
    \label{fig:lbdp1}
  }
\end{figure}

\Cref{fig:lbdp1} shows the results of a calculation such as might form part of a data analysis.
For genealogies induced by the linear birth-death-sampling process, exact expressions for the likelihood are available \citep{Stadler2010}.
Accordingly, we compare the results obtained by integrating \cref{eq:wde1_lbdp,eq:wdiff2_lbdp}, using a sequential Monte Carlo integration scheme \citep{King2016} to these exact results.

\paragraph{SIR model}

\begin{figure}
\begin{knitrout}\small
\definecolor{shadecolor}{rgb}{0.969, 0.969, 0.969}\color{fgcolor}

{\centering \includegraphics[width=1\linewidth]{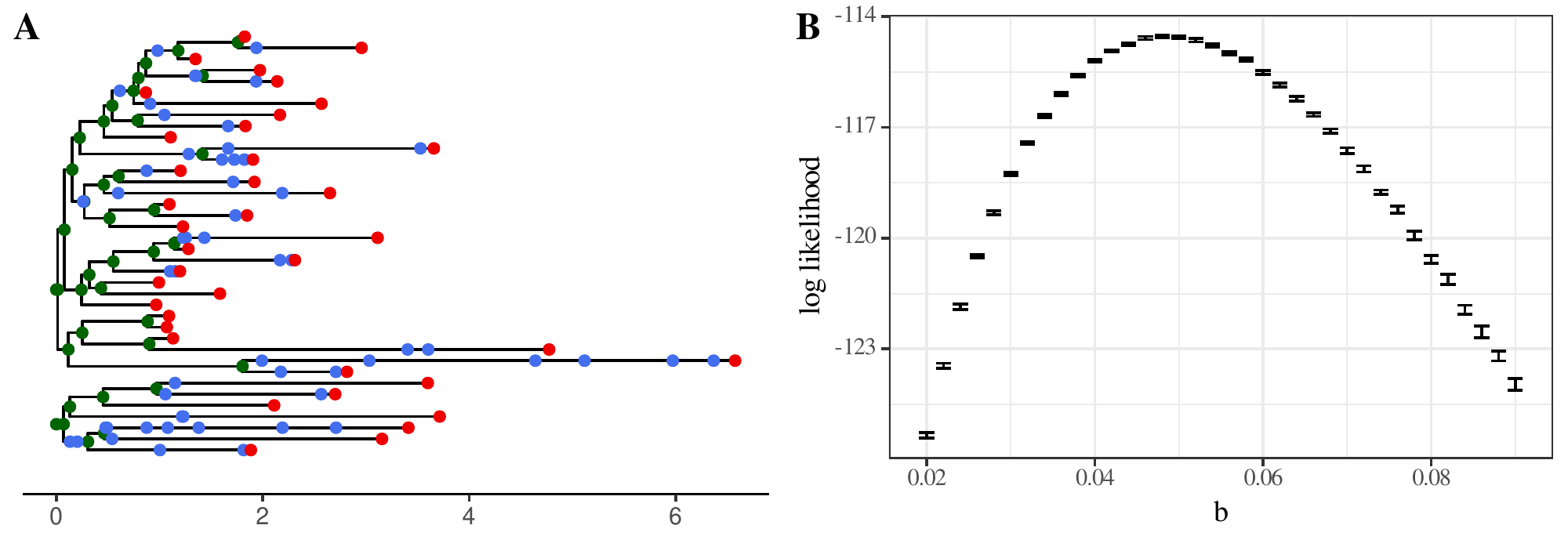} 

}

\end{knitrout}
  \caption{
    \textbf{Computing the likelihood for genealogies induced by the SIR process.}
    \textbf{(A)} A simulated genealogy for $b=0.04$, $\gamma=1$, $\psi=1$.
    The population is of size $100$, with $3$ infectives at time 0.
    \textbf{(B)} Change in the log likelihood as a function of position along a line passing through the true parameter in the $b$-direction.
    The error bars show the estimates obtained using the particle filter (mean ${\pm}2$~s.e.).
    \label{fig:sir1}
  }
\end{figure}

In the case of the SIR model (\cref{sec:examples}), \cref{eq:wde1,eq:wdiff2} become, for $t\notin\Eset(\Vis)$,
\begin{equation*}\label{eq:wde1_sir}
  \begin{split}
    \frac{\partial{w}}{\partial{t}} &= {b\,(s+1)\,(i-1)}\,\left(1-\frac{\binom{\ell}{2}}{\binom{i}{2}}\right)\,w(t,s+1,i-1)\\
    &\quad+{\gamma\,(i+1)}\,w(t,s,i+1)-(b\,s\,i+\gamma\,i+\psi\,i)\,w(t,s,i),
  \end{split}
\end{equation*}
while for $t\in\Eset(\Vis)$,
\begin{equation*}\label{eq:wdiff2_sir}
  w(t,s,i)=
  \begin{cases}
    \frac{b\,(s+1)\,(i-1)}{\mu}\,\frac{1}{\binom{i}{2}}\,w(\prev{t},s+1,i-1), &t\in\Cset(\Vis),\\[13pt]
    \frac{\psi}{\mu}\,w(\prev{t},s,i), &t\in\Qset(\Vis),\\[13pt]
    \frac{\psi}{\mu}\,\left(i-{\ell}\right)\,w(\prev{t},s,i), &t\in\Lset(\Vis).\\
  \end{cases}
\end{equation*}
Again, for simplicity, we have taken all parameters to be constant in time and we understand that $w(t,s,i)=0$ for $i<\ell(t)$.
\Cref{fig:sir1} shows the results of a typical calculation performed using these equations.

\paragraph{SIRS model}

\begin{figure}
\begin{knitrout}\small
\definecolor{shadecolor}{rgb}{0.969, 0.969, 0.969}\color{fgcolor}

{\centering \includegraphics[width=1\linewidth]{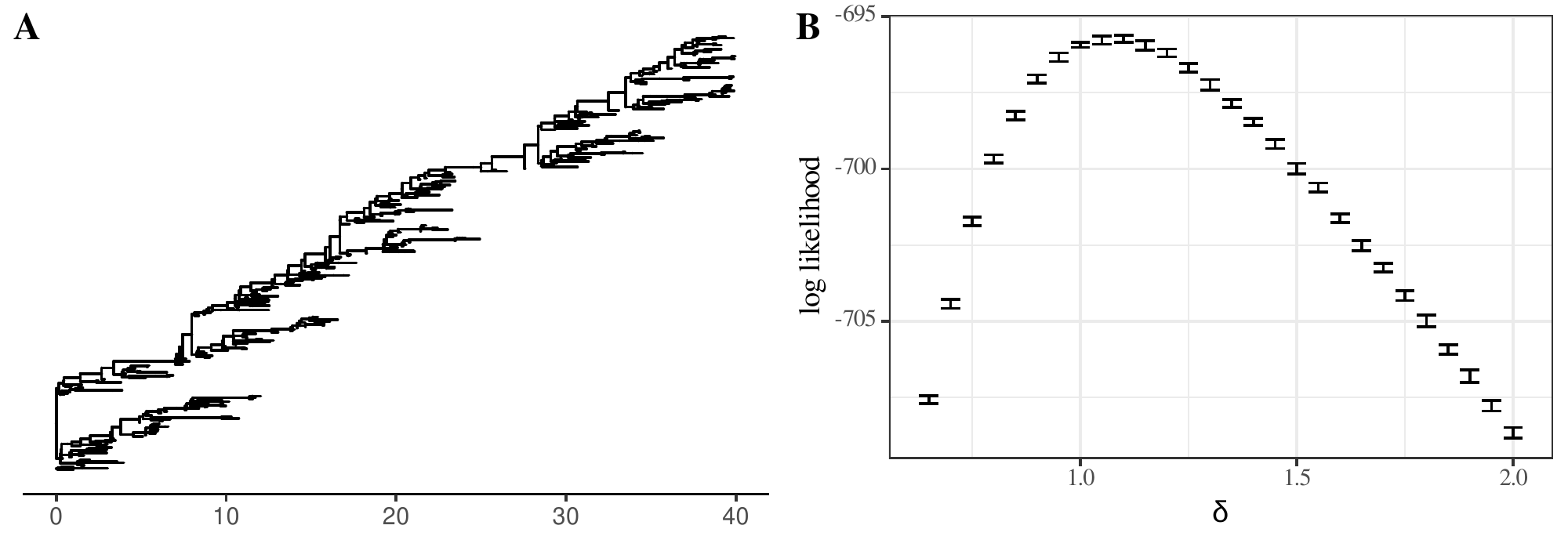} 

}

\end{knitrout}
  \caption{
    \textbf{Computing the likelihood for genealogies induced by the SIRS process.}
    \textbf{(A)} A simulated genealogy for $b=0.04$, $\gamma=2$, $\psi=1$, and $\delta=1$.
    The population is of size $100$, with $3$ infectives at time 0.
    \textbf{(B)} Change in the log likelihood as a function of position along a line passing through the true parameter in the $\delta$-direction.
    The error bars show the estimates obtained using the particle filter (mean ${\pm}2$~s.e.).
    \label{fig:sirs1}
  }
\end{figure}

The case of the SIRS model (\cref{sec:examples}) is similar.
Again, taking all parameters to be constant for simplicity, \cref{eq:wde1,eq:wdiff2} assume the following forms.
For $t\notin\Eset(\Vis)$,
\begin{equation*}\label{eq:wde1_sirs}
  \begin{split}
    \frac{\partial{w}}{\partial{t}} &=
    {b\,(s+1)\,(i-1)}\,\left(1-\frac{\binom{\ell}{2}}{\binom{i}{2}}\right)\,w(t,s+1,i-1,r)
    +{\gamma\,(i+1)}\,w(t,s,i+1,r-1)\\
    &\quad+{\delta\,(r+1)}\,w(t,s-1,i,r+1)
    -(b\,s\,i+\gamma\,i+\delta\,r+\psi\,i)\,w(t,s,i,r).
  \end{split}
\end{equation*}
For $t\in\Eset(\Vis)$, we have
\begin{equation*}\label{eq:wdiff2_sirs}
  w(t,s,i,r)=
  \begin{cases}
    \frac{b\,(s+1)\,(i-1)}{\mu}\,\frac{1}{\binom{i}{2}}\,w(\prev{t},s+1,i-1,r), &t\in\Cset(\Vis),\\[13pt]
    \frac{\psi}{\mu}\,w(\prev{t},s,i,r), &t\in\Qset(\Vis),\\[13pt]
    \frac{\psi}{\mu}\,\left(i-{\ell}\right)\,w(\prev{t},s,i,r), &t\in\Lset(\Vis).\\
  \end{cases}
\end{equation*}
As usual, we have $w(t,s,i,r)=0$ whenever $i<\ell(t)$.
\Cref{fig:sirs1} shows the results of a calculation performed using these equations.

\section{Discussion}
\label{sec:discussion}

Numerical solution of \cref{eq:dmz} is readily achieved, as in \cref{sec:illus}, by applying a Monte Carlo Feynman-Kac approach, i.e., by simulating individual weighted realizations (particles) of the population process between genealogical event times, updating the weights of the particles appropriately at each genealogical event (i.e., according to \cref{eq:wdiff2}).
Importantly, because the genealogical events correspond to real events in the population process, not only the weights, but also the states, of the particles must be adjusted according to \cref{eq:wdiff2}.
Thus, this approach leads to a modified version of the standard sequential Monte Carlo (particle filter) algorithm \citep[e.g.,][]{King2016}.

Indeed, in the special case that the event rates $\alpha$ are time-homogeneous, the sequential Monte Carlo algorithm for evaluating \cref{eq:lik} just described is precisely that proposed by \citet{Vaughan2019}.
Our work here therefore shows how this approach can be extended to a much broader class of models.
A forthcoming paper will extend the class still farther, to encompass simultaneous birth, death, and sampling events, and will describe several alternative algorithms for the numerical solution of \cref{eq:dmz}.

One key feature of the algorithms proposed here is that they enjoy the \emph{plug-and-play property} \citep{He2010}.
An algorithm for inference on partially observed Markov processes is said to be plug-and-play if it operates without ever needing to evaluate the Markov transition probability densities.
Numerical solution of \cref{eq:dmz} requires only that one be able to simulate the population process.
Such methods are attractive inasmuch as they allow consideration of models that are scientifically interesting but mathematically inconvenient by other approaches, since it is typically the case that models can be simulated even when their probability density functions are mathematically intractable.
The fact that plug-and-play methods avoid the need for method-specific approximations also facilitates objective model comparison, since it puts models on an even footing with respect to inference methodology \citep{He2010,King2016}.

As we have mentioned at various points above, it is possible to extend the constructions developed here to more general situations.
In particular, the state space for the population process $\X$ can be taken to be any separable Banach space, so long as the birth, death, sampling, and population size functions, $B$, $D$, $G$, and $I$, remain well defined.
One can also relax the requirements that $B$, $D$, and $G$ have ranges in $\{0,1\}$ and the assumption that these different kinds of events never coincide.
It is necessary to relax these assumptions, for example, if one wishes to entertain models of superspreading \citep{LloydSmith2005a}, or more generally to allow for overdispersion in the latent population process, often an important component of well-fitting models \citep{He2010,Breto2011}.
In accommodating these extensions, the combinatoric arguments of \cref{thm:hardway,thm:easyway} are more intricate, but remain tractable.
We will describe these extensions in a future paper.

The first line of \cref{eq:dmz} resembles the Kolmogorov forward equation (\cref{eq:kfe}) for the population process.
Accordingly, in the absence of sampling, \cref{eq:dmz} preserves the normalization of $w$, i.e., $\sum_x w(t,x,\Vis)=1$.
The remaining lines represent the accumulation of evidence supporting each of the alternative hypotheses $\X_t=x$.
Note that some of the evidence assimilated into $w$ at any time $t$ is derived from data that are only collected \emph{after} time $t$.
Because of this, the partial weights are not measurable in the filtrations induced by the Markov processes of \cref{fig:constellation}.
One must therefore resist the temptation to over-interpret the partial weights $w(t,x,\Vis)$ for $t<t(\Vis)$:
they should be viewed merely as elements of an algorithm that ultimately yields the full likelihood.
Some previous full-information phylodynamic approaches have made mistakes of this kind \citep[e.g.,][]{Rasmussen2011,Leventhal2014}.
The correct interpretation of \cref{eq:dmz} is that a portion of the information in each sample is \emph{referred} to earlier times.
This referral is in some sense the reverse of the evolutionary process whereby information about transmission and recoveries (or births and deaths, or speciations and extinctions) is stored in the genome.
From this point of view, the genealogy itself can be understood as nothing other than a prescription for this information referral.

\section*{Acknowledgments}

The authors gratefully acknowledge useful conversations with Alexandre Bouchard-C\^ot\'e, Simon Frost, Katia Koelle, Vladimir Minin, Mitchell Newberry, David Rasmussen, Jonathan Terhorst, Erik Volz, and two anonymous reviewers.
This work was supported by grants from
the U.S. National Institutes of Health, (Grant \#1R01AI143852 to AAK, \#1U54GM111274 to AAK and ELI)
and a grant from the Interface program, jointly operated by the U.S. National Science Foundation and the National Institutes of Health (Grant \#1761603 to ELI and AAK).
QL was supported by a fellowship from the Michigan Institute for Data Science.

\end{document}